\documentclass{fic-l}

\usepackage{latexsym}
\usepackage{amsthm}
\usepackage[dvips]{graphicx}
\usepackage{epsf}

\newtheorem{theorem}{Theorem}[section]

\newtheorem{proposition}[theorem]{Proposition}

\theoremstyle{definition}
\newtheorem{definition}[theorem]{Definition}
\newtheorem{example}[theorem]{Example}
\newtheorem{corollary}[theorem]{Corollary}

\theoremstyle{remark}
\newtheorem{remark}[theorem]{Remark}

\numberwithin{equation}{section}

\newcommand{\NN}{\ensuremath{\mathbb{N}}}
\newcommand{\ZZ}{\ensuremath{\mathbb{Z}}}
\newcommand{\QQ}{\ensuremath{\mathbb{Q}}}
\newcommand{\RR}{\ensuremath{\mathbb{R}}}
\newcommand{\CC}{\ensuremath{\mathbb{C}}}
\newcommand{\FF}{\ensuremath{\mathbb{F}}}
\newcommand{\GG}{\ensuremath{\mathbb{G}}}
\newcommand{\KK}{\ensuremath{\mathbb{K}}}
\renewcommand{\SS}{\ensuremath{\mathbb{S}}}
\newcommand{\PP}{\ensuremath{\mathbb{P}}}
\newcommand{\OO}{\ensuremath{\mathcal{O}}}
\newcommand{\QQl}{\QQ_{\ell}}
\newcommand{\Fr}{\mathrm{Fr}}
\newcommand{\prim}{\mathrm{prim}}

\newcommand{\ov}{\bar}

\newcommand{\Pic}{\operatorname{Pic}}

\newcommand{\Hdg}{\operatorname{Hdg}}
\newcommand{\End}{\operatorname{End}}
\newcommand{\trace}{\operatorname{Trace}}
\newcommand{\rank}{\operatorname{rank}}

\newcommand{\lam}{\lambda}
\usepackage{epstopdf}

\begin{document}

\title{Modularity of Calabi--Yau varieties: 2011
and beyond}

\author{Noriko YUI}
\address{Department of Mathematics and Statistics\\ Queen's University\\ 
Kingston, Ontario, Canada K7L 3N6 \\ E-mail: yui@mast.queensu.ca}
\thanks{The author was supported in part by NSERC 
Discovery Grant.}


\begin{abstract}
This paper presents the current status on 
modularity of Calabi--Yau varieties since the
last update in 2003. We will focus on Calabi--Yau 
varieties of dimension at most three. 
Here modularity refers to at least two different
types: arithmetic modularity and geometric modularity.
These will include: (1) the modularity (automorphy) 
of Galois representations of Calabi--Yau varieties (or motives) 
defined over $\QQ$ or number fields, (2) the modularity of solutions of 
Picard--Fuchs differential equations of families of Calabi--Yau varieties,
and mirror maps (mirror moonshine), (3) the modularity of
generating functions of invariants counting certain quantities
on Calabi--Yau varieties, and (4) the modularity of moduli for 
families of Calabi--Yau varieties.  The topic (4) is commonly
known as geometric modularity.  

Discussions in this paper are centered around  
arithmetic modularity, namely on (1), and (2), with a brief excursion 
to (3).
\end{abstract}

\maketitle

\section*{Introduction}
These are notes of my introductory lectures at the Fields workshop
on ``Arithmetic and Geometry of K3 surfaces and
Calabi--Yau threefolds'', at the Fields
Institute from August 16 to August 25, 2011.
My goal is two-fold: 

\begin{itemize}
\item Present an update on the recent developments on 
the various kinds of modularity associated to Calabi--Yau
varieties. Here ``recent'' developments means various developments 
since my article \cite{Yui2003} published in 2003.
\item Formulate conjectures and identify future problems on
modularity and related topics.
\end{itemize}

My hope is to motivate young (as well as mature) researchers 
to work in this fascinating area at the interface of 
arithmetic, geometry and physics around Calabi--Yau 
varieties.

\subsection{Brief history since 2003}

The results and discoveries in the last ten years on 
Calabi--Yau varieties, which will be touched upon in my
lectures, are listed below.  

\begin{itemize}
\item The modularity of the $2$-dimensional
Galois representations associated to  
Calabi--Yau varieties defined over $\QQ$. 
\item The modularity of highly reducible Galois representations 
associated to Calabi--Yau threefolds over $\QQ$.
\item The automorphy of higher dimensional
Galois representations arising from CM type Calabi--Yau
varieties (automorphic induction).
\item  Appearance of various types of modular forms 
in Mirror Symmetry as generating functions counting
some mathematical/physical quantities.
\item  Modularity of families of Calabi--Yau varieties
(solutions of Picard--Fuchs differential equations, monodromy 
groups, mirror maps).
\item Moduli spaces of Calabi--Yau families, and higher
dimensional modular forms (e.g., Siegel modular forms).
\end{itemize}

\subsection{Plan of lectures}

Obviously, due to time constraints, I will not be able to cover all
these topics in my two introductory lectures. Thus, 
for my lectures, I plan to focus on recent results on the
first three (somewhat intertwined) items listed above, and 
with possibly very brief interludes to
the rest if time permits; 

This article includes my two lectures delivered at
the workshop, as well as some subjects/topics which I was not able to 
cover in my lectures.  However, I must emphasize that this note will not
touch upon geometric modularity. 

\subsection{Disclaimer}
Here I will make the disclaimer that the topics listed
above are by no means exhaustive; it may be the case
that I forget to mention some results, or not
give proper attributions.  I apologize for these oversights.

\subsection{Calabi--Yau varieties: Definition}

\begin{definition} Let $X$ be a smooth projective
variety of dimension $d$ defined over $\CC$. We say
that $X$ is a {\it Calabi--Yau} variety if
\begin{itemize}
\item  $H^i(X,\mathcal{O}_X)=0$ for every $i,\,0<i<d$. 
\item  The canonical bundle $K_X$ is trivial.
\end{itemize}

We introduce the {\it Hodge numbers} of $X$:
$$h^{i,j}(X):=\mathrm{dim}_{\CC} H^j(X,\Omega_X^i),\,
0\leq i, j\leq d.$$

Then we may characterize a Calabi--Yau variety of dimension $d$
in terms of its Hodge numbers.

A smooth projective variety $X$ of
dimension $d$ over $\CC$ is called a {\it Calabi--Yau}
variety if
\begin{itemize}
\item  $h^{i,0}(X)=0$ for every $i, \,0 < i <d.$
\item  $K_X\simeq \mathcal{O}_X$, so that
the geometric genus of $X$, $$p_g(X):=h^{0,d}(X)
=\mathrm{dim}_{\CC} H^0(X, K_X)=
\mathrm{dim}_{\CC}H^0(X,\mathcal{O}_X)=1.$$
\end{itemize}
\end{definition}

\noindent{\bf Numerical characters of Calabi--Yau
varieties of dimension $d$}

\begin{itemize}
\item  Betti numbers: For $i,\, 0\leq i\leq 2d$,
the $i$-th Betti number of $X$ is defined by
$$B_i(X):=\mathrm{dim}_{\CC} H^i(X,\CC).$$
There is Poincar\'e duality for $H^i(X,\CC)$;
that is,
$$H^i(X,\CC)\times H^{2d-i}(X,\CC)\to
\CC\quad\mbox{for every $i,\,0\leq i \leq d$}$$
is a perfect pairing. This implies that
$$B_i(X)=B_{2d-i}(X),\mbox{for $i,\,0\leq i\leq d$}$$
\item  Hodge numbers: They are defined in Definition 0.1 above. There is
the symmetry of Hodge numbers: For $0\leq i,j\leq d$,
$$h^{i,j}(X)=h^{j,i}(X)\quad\mbox{by complex conjugation},$$ and
$$h^{i,j}(X)=h^{d-i,d-j}(X)\quad\mbox{by Serre duality}.$$
\item  There is a relation among Betti numbers
and Hodge numbers, as a consequence of the Hodge decomposition:
$$B_k(X)=\sum_{i+j=k} h^{i,j}(X).$$
\item  The Euler characteristic of $X$ is defined by
$$E(X):=\sum_{k=0}^{2d} (-1)^k B_k(X).$$
\end{itemize}

\begin{example} 
(a) Let $d=1$. The first condition is vacuous. The
second condition says that $p_g(X)=1$. So 
dimension $1$ Calabi--Yau varieties are elliptic curves.
The Hodge diamond of elliptic curves is rather simple.

$$
\begin{array}{c@{}c@{}c@{}c@{}c@{}cl}
&&  1 && \qquad\qquad &B_0(X)=1 \\
1\phantom{123} &&&&\phantom{123} 1 \qquad\qquad &B_1(X)=2 \\
&& 1 &&  \qquad\qquad &B_2(X)=1 \\
\end{array}
$$
The Euler characteristic of $X$ is given by
$$E(X)=B_0(X)-B_1(X)+B_2(X)=0.$$

(b) Let $d=2$. The first condition is $h^{1,0}(X)=0$
and the second condition says that $h^{0,2}(X)=p_g(X)=1$. 
So dimension $2$ Calabi--Yau varieties
are K3 surfaces.
The Hodge diamond of K3 surfaces is of the
form:

$$\begin{array}{c@{}c@{}c@{}c@{}c@{}c@{}c@{}c@{}cl}
&& & 1 & & &\qquad\qquad &B_0(X)=1 \\
& 0\phantom{123} & &&&\phantom{123} 0 & \qquad\qquad &B_1(X)=0 \\
1\phantom{123} & &&20 && &\phantom{123} 1 \qquad\qquad &\,\,\, B_2(X)=20 \\
& 0\phantom{123} &&&&\phantom{123} 0 & \qquad\qquad & B_3(X)=0 \\
&& & 1 && & \qquad\qquad &B_4(X)=1 \\
\end{array}
$$
The Euler characteristic of $X$ is 
given by
$$E(X)=\sum_{k=0}^4 (-1)^kB_k(X)=1+22+1=24.$$

(c) Let $d=3$. The first condition says that
$h^{1,0}(X)=h^{2,0}(X)=0$ and the second condition
implies that $h^{0,3}(X)=p_g(X)=1$.  
The Hodge diamond of $X$ is given as follows:

$$
\begin{array}{c@{}c@{}c@{}c@{}c@{}c@{}c@{}cl}
  &   &       &          1 &            &   &   & \qquad\qquad &
 B_0(X) = 1\\
  &   & 0          &            &    0 &   &   & & B_1(X) = 0\\

  & 0 &            & h^{1,1}(X) &   & 0 &   & & B_2(X) = h^{1,1}(X)\\
1\phantom{123} &   & h^{2,1}(X) &            & h^{1,2}(X) &   & \phantom{123}1 & & B_3(X) = 2(1+h^{2,1}(X))\\
  & 0 &            & h^{2,2}(X) &            & 0 &   & & B_4(X) = h^{2,2}(X) = h^{1,1}(X)\\
  &   & 0          &            &          0 &   &   & & B_5(X) = 0\\
  &   &            &          1 &       &   &   & & B_6(X) = 1
\end{array}
$$

The Euler characteristic is given by
$$E(X)=\sum_{k=0}^6 (-1)^k B_k(X)=2(h^{1,1}(X)-h^{2,1}(X)).$$
It is not known if there exist absolute constants that
bound $h^{1,1}(X)$, $h^{2,1}(X)$, and hence $|E(X)|$.  The
currently known bound for $|E(X)|$ is $960$.

(d) Here are some typical examples of families of Calabi--Yau
varieties defined by hypersurfaces.

$$\begin{array}{c||c||c}  \hline
d & \mbox{$CY$ variety of dim $d$}& \mbox{$CY$ varieties of dim $d$} \\ \hline
1 & X_0^3+X_1^3+X_2^3+3\lam X_0X_1X_2 & y^2=f_3(x) \\ \hline
2 & X_0^4+X_1^4+X_2^4+X_3^4+4\lam X_0X_1X_2X_3 & z^2=f_6(x,y)\\ 
\hline
3 & X_0^5+X_1^5+X_2^5+X_3^5+X_4^5+5\lam X_0X_1X_2X_3X_4 &
w^2=f_8(x,y,z) \\ \hline
\end{array}$$
\medskip

The equations in the second column are generic polynomials in projective 
coordinates;
while those in the third column are in affine coordinates and
the $f_i$ are smooth polynomials of degree $i$ in affine coordinates.  

We have a vast source of examples of Calabi--Yau threefolds via toric 
construction and other methods. The upper record for the absolute value
of the Euler characteristic of all these
Calabi--Yau threefolds is $960$, though there is neither reason
nor explanation for this phenomenon.
\end{example}

\section{The modularity of Galois representations of
Calabi--Yau varieties (or motives) over $\QQ$.}

We now consider smooth projective varieties 
defined over $\QQ$ (say, by hypersurfaces, or
complete intersections).
We say that $X/\QQ$ is a Calabi--Yau variety of dimension $d$
over $\QQ$, if $X\otimes_{\QQ}\CC$ is a Calabi--Yau variety
of dimension $d$. 
A Calabi--Yau variety $X$ over $\QQ$ has a model defined 
over $\ZZ[\frac{1}{m}]$ (with some $m\in\NN$), and this allows 
us to consider its reduction modulo primes.  Pick a prime
$p$ such that $(p,m)=1$, and define the reduction of $X$ modulo $p$, denoted by
$X\mod p$.  We say that $p$ is a {\it good} prime if $X\mod p$
is smooth over $\ov{\FF}_p$, otherwise $p$ is {\it bad}. 
For a good prime $p$, let $\Fr_p$ denote
the Frobenius morphism on $X$ induced from 
the $p$-th power map $x\mapsto x^p$.

Let $\ell$ be a prime different from $p$. Then,
for each $i, \, 0\leq i\leq 2d$, $\Fr_p$ induces
an endomorphism $\Fr_p^*$ on the $i$-th $\ell$-adic \'etale 
cohomology group $H^i_{et}(\ov{X}_p,\QQl)$,
where $\ov{X}_p:=X\otimes_{\FF_p}\ov{\FF}_p$.
Grothendieck's specialization theorem gives
an isomorphism
$H^i_{et}(\ov{X}_p,\QQl)\cong H^i_{et}(\ov{X},\QQl)$,
where $\ov{X}:=X\otimes_{\QQ}\ov{\QQ}$. 
Then the comparison theorem gives 
$H^i_{et}(\ov{X},\QQl)\otimes_{\QQl}\CC\cong
H^i(X\otimes_{\QQ} \CC, \CC)$ so that $\dim_{\QQl} H^i_{et}(\ov{X},\QQl)=B_i(X).$
There is Poincar\'e duality for $H^i_{et}(\ov{X},\QQl)$, 
that is,
$$H^i_{et}(\ov{X},\QQl)\times H^{2d-i}_{et}(\ov{X},\QQl)\to
\QQl\quad\mbox{for every $i,\,0\leq i\leq 2d$}$$
is a perfect pairing. 
Let
$$P_p^i(T):=\mbox{det}(1-\Fr_p^*\,T\,|\, 
H^i_{et}(\ov{X}, \QQl))$$
be the reciprocal characteristic polynomial of $\Fr_p^*$.  (Here
$T$ is an indeterminate.)  Then the Weil Conjecture
(Theorem) asserts that

\begin{itemize}
\item $P_p^i(T)\in 1+T\ZZ[T]$. Moreover,
$P_p^i(T)$ does not depend on the choice of $\ell$. 
\item $P_p^i(T)$ has degree $B_i(X)$.
\item $P_p^{2d-i}(T)=\pm P_p^i(p^{d-\frac{i}{2}}T)$ for every 
$i,\, 0\leq i\leq d$.
\item If we write
$$P_p^i(T)=\prod_{j=1}^{B_i} (1-\alpha_{ij}\,T)\in\ov{\QQ}[T]$$
then $\alpha_{ij}$ are algebraic integers with $|\alpha_{ij}|=p^{i/2}$ for
every $i,\,0\leq i\leq 2d$.
\end{itemize}

Now we will bring  in the absolute Galois group
$G_{\QQ}:=\mbox{Gal}(\ov{\QQ}/\QQ)$.  There is
a continuous system of $\ell$-adic Galois representations
$$\rho^i_{X,\ell}: G_{\QQ}\to GL(H^i_{et}(\ov{X},\QQl))$$
sending the (geometric) Frobenius
$\Fr_p^{-1}$ to $\rho^i(\Fr_p^{-1})$. The 
(geometric) Frobenius $\rho^i(\Fr_p^{-1})$ has the same
action as the Frobenius morphism $\Fr_p^*$
on the \'etale cohomology $H^i_{et}(\ov{Z},\QQl)$.
We define its $L$-series 
$L(\rho^i_{X,\ell}, s):=L(H^i_{et}(\ov{X},\QQl), s)$
for each $i,\, 0\leq i\leq 2d$,
where $s$ is a complex variable. 

We will now define the $L$-series of $X$. 

\begin{definition} The $i$-th (cohomological) $L$-series
of $X$ is defined by the Euler product
$$L_i(X,s):=L(H^i_{et}(\ov{X},\QQl),s)$$
$$:=(*)\prod_{p: p\neq\ell}
P_p^i(p^{-s})^{-1}
\times (\mbox{factor corresponding to $\ell=p$})$$
where the product is taken over all good primes
different from $\ell$ and $(*)$ corresponds to factors
of bad primes. 
For $\ell=p$, we may choose another good prime $\ell\neq p$,
or we can use some $p$-adic cohomology
(e.g., crystalline cohomology) to define the factor.
\medskip

For $i=d$, we write simply $L(X,s)$ for $L_d(X,s)$
if there is no danger of ambiguity.
\end{definition}

\begin{remark}
We may define (for a good prime $p$) the zeta-function 
$\zeta(X_p,T)$, of a Calabi--Yau variety $X_p$
defined over $\FF_p$ by counting the number of
rational points on all extensions of $\FF_p$:
$$\zeta(X_p,T):=\mathrm{exp}\left(\sum_{n=1}^{\infty}
\frac{\# X_p(\FF_{p^n})}{n}T^n\right)\in\QQ(T).$$
Then by Weil's conjecture, it has the form: 
$$\zeta(X_p,T)=\frac{P_1(T)P_3(T)\cdots P_{2d-1}(T)}{P_0(T)P_2(T)\cdots P_{2d}(T)},$$ 
where we put (to ease the notation) $P_i(T)=P_p^i(T)$ for $i=1,\ldots, 2d$.
\end{remark}
\medskip

Let $X$ be a Calabi--Yau variety of dimension $d$
defined over $\QQ$. For a good prime $p$, $\Fr_p$ acts on $X\mod p$,
and it will induce a morphism $\Fr_p^*$ on
$H^i_{et}(\ov{X}_p,\QQl)\simeq H^i_{et}(\ov{X},\QQl)$. 
Define the trace $t_i(p)$ by 
$$t_i(p):=\trace(\Fr_p^*\,|\,H^i_{et}(\ov{X},\QQl)),
\quad \mbox{for $i,\,0\leq i\leq 2d$}.$$
Then $t_i(p)\in\ZZ$ for every $i, \, 0\leq i\leq 2d$.
The Lefschetz fixed point formula gives a relation
between the number of $\FF_p$-rational points on $X$
and traces:
$$\#X(\FF_p)=\sum_{i=0}^{2d} (-1)^i t_i(p).$$

\begin{example}
(a) Let $d=1$ and let $E$ be an elliptic curve defined
over $\QQ$. For a good prime $p$, 
$$\#E(\FF_p)=\sum_{i=0}^2 (-1)^i t_i(p)=t_0(p)-t_1(p)+t_2(p)
=1+p-t_1(p),$$
where
$$|t_1(p)|\leq 2p^{1/2}.$$
Then we have
$$
P_p^0(T)=1-T,\, P_p^2(T)=1-pT,\, P_p^1(T)=1-t_1(p)T+pT^2.$$
The $L$-series of $E$ is then given by
$$L(E,s)=(*)\prod_{p:good} P^1_p(p^{-s})^{-1}$$
$$=(*)\prod_{p:good}\frac{1}{1-t_1(p)p^{-s}+p^{1-2s}}.$$
Expanding out, we may write
$$L(E,s)=\sum_{n=1}^{\infty}\frac{a(n)}{n^s}\quad
\mbox{with $a_1=1$ and $a(n)\in\ZZ$.}$$
So $$a(p)=t_1(p)\quad\mbox{for every good prime $p$.}$$

(b) Let $d=2$, and let $X$ be a K3 surface defined
over $\QQ$. Let $NS(X)$ denote the N\'eron--Severi
group of $X$ generated by algebraic cycles. It is
a free finitely generated abelian group, and
$$NS(X)=H^2(X,\ZZ)\cap H^{1,1}(X)$$
so that the rank of $NS(X)$ (called the Picard
number of $X$), denoted by $\rho(X)$, 
is bounded above by $20$. Let $T(X)$ be the 
orthogonal complement of $NS(X)$ in $H^2(X,\ZZ)$ with
respect to the intersection pairing. We call 
$T(X)$ the group of transcendental cycles of $X$.
We have the decomposition
$$H^2(X,\ZZ)\otimes\QQl=(NS(X)\otimes\QQl)\oplus (T(X)\otimes\QQl).$$
and this will enable us decompose the $L$-series of $X$
as follows:
$$L(X,s)=L(H^2_{et}(\ov{X},\QQl),s)
=L(NS(X)\otimes\QQl,s)\times L(T(X)\otimes\QQl,s).$$

We know that for a good prime $p$,
$$P_p^0(T)=1-T, P_p^4(T)=1-p^2T,\, P_p^1(T)=P_p^3(T)=1$$
and if
$$P_p^2(T)=\prod_{j=1}^{22}(1-\alpha_j\,T)\in\ov{\QQ}[T]$$
then $$|\alpha_j|=p.$$

The $L$-series of $NS(X)$ is more or less understood
by Tate's conjecture. The validity of the Tate conjecture for K3 surfaces
in characteristic zero has been established (see
Tate \cite{T94}). In fact, if we know that 
all the algebraic cycles generating $NS(X)$ are defined over some
finite extension $\KK$ of $\QQ$ of degree $r$,
then $\rho^2(\Fr_{p^r})$ acts on $NS(X)\otimes \QQl$ by
multiplication by $p^r$ so that the $L$-series
may be expressed as
$$L(NS(X_{\KK})\otimes\QQl,s)=\zeta_{\KK}(s-1)^{\rho(X)}$$
where $\zeta_{\KK}(s)$ denotes the Dedekind zeta-function
of $\KK$.  

Therefore, the remaining task is to determine
the $L$-series $L(T(X)\otimes\QQl,s)$ arising from the transcendental cycles $T(X)$,
and we call this the motivic $L$-series of $X$.

(c) Let $d=3$, and let $X$ be a Calabi--Yau threefold
defined over $\QQ$.  For a good prime $p$, 
$$\#X(\FF_p)=\sum_{i=0}^6 (-1)^i t_i(p)
=t_0(p)-t_1(p)+t_2(p)-t_3(p)+t_4(p)-t_5(p)+t_6(p)$$
$$=1+p^3+(1+p)t_2(p)-t_3(p).$$
Therefore,
$$t_3(p)=1+p^3+(1+p)t_2(p)-\#X(\FF_p).$$
We have $$
t_0(p)=1,\,t_6(p)=p^3,\, t_1(p)=t_5(p)=0,$$
$$|t_2(p)|\leq ph^{1,1}(X),\, t_4(p)=pt_2(p),$$
and
$$|t_3(p)|\leq B_3 \,p^{3/2}.$$
Hence
$$P_p(T)=1-T,\,P_p^6(T)=1-p^3T,\, P_p^1(T)=P_p^5(T)=1,$$
$$P_p^4(T)=P_2^2(pT)\quad\mbox{and}\quad
P_p^3(T)\in\ZZ[T]\quad\mbox{with degree $B_3(X)$}.$$
Then the $L$-series of $X$ is given by
$$L(X,s)=L(H^3_{et}(\ov{X},\QQl),s)=(*)\prod_{p\,\, good} P_p^3(p^{-s})^{-1}$$
where $(*)$ is the factor corresponding to bad primes.
\end{example}

We will now make a definition of what it means for
a Calabi--Yau variety $X$ defined over $\QQ$ to
be modular (or automorphic). This is a concrete
realization of the conjectures known as the Langlands
Philosophy. 

In the Appendix, we will briefly recall the definitions of various 
types of modular forms. 
\vskip 0.5cm

Now we will recall the Fontaine--Mazur conjectures, or rather
some variant concentrated on Calabi--Yau varieties over $\QQ$.

\begin{definition}
Let $X$ be a Calabi--Yau variety of dimension $d\leq 3$ 
defined over $\QQ$. Let $L(X,s)$ be its $L$-series.
We say that $X$ is modular if there is a set of modular
forms (or automorphic forms) such that
$L(X,s)$ coincides with the $L$-series
associated to modular forms (automorphic forms),
up to a finite number of Euler factors.
\end{definition}

\begin{remark}
In fact, when the Galois representation arising from $X$
is reducible, then we will consider its irreducible factors
and match their $L$-series with the $L$-series of
modular forms. This is the so-called motivic modularity.
\end{remark}
\medskip

Now we define Calabi--Yau varieties $X$ of {\it CM-type}. 
This involves a polarized rational Hodge structure 
$h$ on the primitive cohomology
$H^d_{\prim}(X,\QQ)$, where $d=\mathrm{dim}(X)$.

Let $\SS:=R_{\CC/\RR} \GG_m$ be the real algebraic
group obtained from $\GG_m$ by restriction of scalars
from $\CC$ to $\RR$.
For $d=1$, $H^1_{\prim}(X,\QQ)=H^1(X,\QQ)$,
and polarized rational Hodge structures on $X$ are 
simple.  For $d=2$, 
The Hodge group of a polarized rational Hodge structure
$h: \SS\to GL(H^d_{\prim}(X,\QQ)\otimes\RR)$ 
is the smallest algebraic group 
of $GL(H^2_{\prim}(X,\QQ)\otimes\RR)$
defined over $\QQ$ such that the real points $\Hdg(\RR)$
contain $h(U^1)$ where $U^1:=\{\,z\in\CC^*\,|\, z\overline{z}=1\}$.
For details about Hodge groups, see Deligne
\cite{D900}, and Zarhin \cite{Za83}.
For $d=3$, the Hodge structure on $H^3_{\prim}(X,\QQ)=H^3(X,\QQ)$ is not simple for
all Calabi--Yau threefolds $X$. (For instance, for the quintic
threefold in the Dwork pencil, the Hodge structure would split into
$4$-dimensional Hodge substructures, using the action of the
$(\ZZ/5\ZZ)^3$.)
\begin{definition}
A Calabi--Yau variety $X$ of dimension $d\leq 3$ is said to be of 
CM {\it type} if the Hodge group $\Hdg(X)$ associated to
a rational Hodge structure of weight $k$ 
of $H^d_{\prim}(X,\QQ)$ is commutative.  
That is, the Hodge group $\Hdg(X)_{\CC}$ is isomorphic
to a copy of $\GG_m\simeq \CC^*$.
\end{definition}

Hodge groups are very hard to compute in practice.
Here are algebraic characterizations of CM type
Calabi--Yau varieties of dimension $d$.

\begin{proposition}

$\bullet$ $d=1.$ An elliptic curve $E$ over $\QQ$ is 
of {\mbox{CM}} type if and only if $\End(E)\otimes \QQ$ 
is an imaginary quadratic field over $\QQ$.

$\bullet$ $d=2.$ A {\mbox{K3}} surface $X$ over $\QQ$ is of {\mbox{CM}} type
if 
$\End_{\Hdg}(T(X))\otimes\QQ$ is 
a {\mbox{CM}} field over $\QQ$ of degree equal to $\mathrm{rank}\, T(X)$.

$\bullet$ $d=3.$ Let $X$ be a Calabi--Yau threefold
over $\QQ$.  A Calabi--Yau threefold $X$ over $\QQ$ is of
{\mbox{CM}} type if and only if 
$\End_{\Hdg}(X)\otimes\QQ$
is a {\mbox{CM}} field over $\QQ$ of degree $2(1+h^{2,1}(X))$,
if and only if the Weil and Griffiths intermediate Jacobians of $X$ are 
of {\mbox{CM}} type. 
\end{proposition}

For $d=1$, this is a classical result. For $d=2$, a best reference
might be Zarhin \cite{Za83}, and for $d=3$, see Borcea \cite{Bo92}. 

Later we will construct Calabi--Yau varieties
of dimension $2$ and $3$ which are of CM type.

\section{Results on modularity of Galois representations}

\subsection{Two-dimensional Galois representations arising
from Calabi--Yau varieties over $\QQ$}
\medskip

We will focus on $2$-dimensional Galois representations
arising from Calabi--Yau varieties over $\QQ$.

First, for dimension $1$ Calabi--Yau varieties over $\QQ$,
we have the celebrated theorem of Wiles et al.

\begin{theorem} ($d=1$) Every elliptic curve $E$
defined over $\QQ$ is modular. More concretely, let
$E$ be an elliptic curve over $\QQ$ with conductor $N$.
Then there exists a Hecke eigen newform $f$ of weight $2=1+d$
on the congruence subgroup $\Gamma_0(N)$ such that
$$L(E,s)=L(f,s).$$
That is, if we write $f(q)=\sum_{m=1}^{\infty} a_f(m)q^m$
with $q=e^{2\pi iz}$ and normalized by $a_f(1)=1$, then
$a(n)=a_f(n)$ for every $n$.
\end{theorem}

For dimension $2$ Calabi--Yau varieties, namely,
K3 surfaces, over $\QQ$, there are naturally associated 
Galois representations. In particular, for a special class of K3 surfaces,
the associated Galois representations are $2$-dimensional.  
Let $X$ be a K3 surface
defined over $\QQ$ with Picard number $\rho(X_{\ov{\QQ}})=20$.
Such K3 surfaces are called {\it singular} K3 surfaces. Then
the group (or lattice) $T(X)$ of transcendental cycles on $X$
is of rank $2$, and it gives rise to a $2$-dimensional
Galois representations. Livn\'e \cite{L95} has established
the motivic modularity of $X$, that is, the modularity of $T(X)$. 

\begin{theorem} ($d=2$) Let $X$ be a singular K3 surface
defined over $\QQ$. Then $T(X)$ is modular, that is,
there is a modular form $f$ of weight $3=1+d$ on some
$\Gamma_1(N)$ or $\Gamma_0(N)$ with a character $\varepsilon$ 
such that
$$L(T(X)\otimes\QQl,s)=L(f, s).$$
\end{theorem}

\begin{remark}
A representation theoretic formulation of the above theorem
is given as follows. Let $\pi$ be the compatible
family of $2$-dimensional $\ell$-adic Galois representations
associated to $T(X)$ and let $L(\pi, s)$ be its
$L$-series. Then there exists a unique, up to isomorphism,
modular form of weight $3$, level=conductor of $\pi$,
and Dirichlet character $\varepsilon(p)=\left(\frac{-d}{p}\right)$ such
that $$L(\pi,s)=L(f,s).$$  Here $d=|\mbox{disc}\,\, NS(X)|$.
\end{remark}

For dimension $3$ Calabi--Yau varieties over $\QQ$,
we will focus on rigid Calabi--Yau threefolds. A
Calabi--Yau threefold $X$ is said to be {\it rigid}
if $h^{2,1}(X)=0$ so that $B_3(X)=2$. This gives
rise to a $2$-dimensional Galois representation.

\begin{theorem}($d=3$)
Every rigid Calabi--Yau threefold $X$ over $\QQ$ is
modular. That is, there exists a cusp 
form $f$ of weight $4=1+d$ on some $\Gamma_0(N)$
such that
$$L(X,s)=L(f,s).$$
\end{theorem} 

This theorem has been established by Gouv\^ea and Yui
\cite{GY2011}, and independently by Dieulefait \cite{D2010}). 
Their proof relies heavily on the recent results on the modularity
of Serre's conjectures about $2$-dimensional
residual Galois representations by Khare--Wintenberger \cite{KW2009}
and Kisin \cite{Ki2009}.

A list of rigid Calabi--Yau threefolds
over $\QQ$ can be found in the monograph of Meyer \cite{M05}. 

\subsection{Modularity of higher dimensional
Galois representations arising from K3 surfaces 
over $\QQ$}

We will first consider K3 surfaces $X$ with $T(X)$
of rank $\geq 3$.  

The first result is for K3 surfaces with transcendental rank $3$.

\begin{theorem} {\sl Let $X$ be a {\mbox{K3}} surface over $\QQ$ with
Picard number $19$.  
Then $X$ has a Shioda--Inose structure,
that is, $X$ has an involution $\iota$ such that $X/\iota$ is
birational to a Kummer surface $Y$ over $\CC$. 

Suppose that the Kummer surface
is given by the product $E\times E$ of a non-{\mbox{CM}} elliptic curve
$E$ over $\QQ$. Then the Shioda--Inose structure induces
an isomorphism of integral Hodge structures on the
transcendental lattices, so, $X$ and ${\mbox{Km}}(E\times E)$ have the
same $\QQ$-Hodge structure.  In this case, the $2$-dimensional
Galois representation $\rho_E$ associated to $E$ induces
the $3$-dimensional Galois representation ${\mbox{Sym}}^2 \rho_E$
on $T(X)$ over some number field. 
Consequently, $T(X)$ is {\bf potentially modular} in the sense
that the $L$-series of $T(X)$ is determined over some
number field $K$ by the symmetric square of a modular form 
$g$ of weight $2$ associated to $E$, and  
over $K$, $$L(T(X)\otimes\QQl,s)=L({\mbox{Sym}}^2(g), s).$$}
\end{theorem}

\begin{remark}
In the above theorem, we are not able to obtain the
modularity results over $\QQ$. 
Since the K3 surface $S$ is defined over $\QQ$, the
representation on $T(X)$ is also defined over $\QQ$,
but the isomorphism to $\mbox{Sym}^2 \rho_E$ may not
be. Thus, we only have the potential modularity of
$X$. 
\end{remark}
\medskip
 
Can we say anything about the modularity of
K3 surfaces with arbitrary large transcendental rank?
The Galois representations associated to these
K3 surfaces have large dimensions.  The only result
along this line is due to Livn\'e-Sch\"utt-Yui
\cite{LSY2010}.

Consider a K3 surface $X$ with non-symplectic
automorphism.  Let $\omega_X$ denote a holomorphic
$2$-form on $X$, fixed once and for all. Then
$H^{2,0}(X)\simeq\CC\omega_X$. 
Let $\sigma\in\mathrm{Aut}(X)$. Then $\sigma$ induces
a map 
$$\sigma^*: H^{2,0}(X)\to H^{2,0}(X)\quad \omega_X\to 
\alpha\omega_X,\,\alpha\in\CC^*.$$
We say that $\sigma^*$ is non-symplectic if $\alpha\neq 1$.
Let $$H_X:=\mathrm{Ker}(\mathrm{Aut}(X)\to O(NS(X)).$$
Then $H_X$ is a finite cyclic group, and in fact, can
be identified with the group of roots of unity $\mu_k$ for some 
$k\in\NN$. Assume that $\mathrm{det}(T(X))=\pm 1$ (that
is, assume that $T(X)$ is unimodular.)  Then we have the
following possibilities for the values of $k$.

$\bullet$ $k\leq 66$ (Nikulin \cite{N80})

$\bullet$ $k$ is a divisor of $66, 44, 42, 36, 28, 12$
(Kondo \cite{K92})

$\bullet$ If $\mathrm{rank}(T(X))=\phi(k)$ (where
$\phi$ is the Euler function), then
$k=66, 44, 42, 36, 28, 12$. Furthermore, there is
a unique K3 surface $X$ with given $k$. (Kondo \cite{K92})
These results were first announced by Vorontsov \cite{V83}, and
proofs were given later by Nikulin \cite{N80} and Kondo \cite{K92}.

We tabulate these six K3 surfaces.

$$\begin{array}{|c|| c|c|c|}  \hline
k & NS(X) & T(X) & \mathrm{rank}(T(X) \\ \hline
12 & U\oplus (-E_8)^2 & U^2 & 4 \\
28 & U\oplus (-E_8) & U^2\oplus (-E_8) & 12 \\
36 & U\oplus (-E_8) & U^2\oplus (-E_8) & 12 \\
42 & U\oplus (-E_8) & U^2\oplus (-E_8) & 12 \\
44 & U & U^2\oplus (-E_8)^2 & 20 \\
66 & U & U^2\oplus (-E_8)^2 & 20 \\ \hline 
\end{array}
$$
\medskip

Here $(-E_8)$ denotes the negative definite even unimodular
lattice of rank $8$.

Now we ought to realize these K3 surfaces over $\QQ$.
Here are explicit equations thanks to Kondo \cite{K92}.

$$\begin{array}{|c|| c| c|}  \hline
k & X & \sigma \\ \hline
12 & y^2=x^3+t^5(t^2+1) & (x,y,t)\mapsto (\zeta^2_{12}x,
\zeta_{12}^3y,-t) \\  
28 & y^2=x^3+x+t^7 & (x,y,t)\mapsto (-x, \zeta_{28}^7y,\zeta_{28}^2t) \\
36 & y^2=x^3-t^5(t^6-1) & (x,y,t)\mapsto (\zeta^2_{36}x,\zeta_36^3y,\zeta_{36}^{30}t) \\
42 & y^2=x^3+t^5(t^7-1) & (x,y,t)\mapsto (\zeta^2_{42}x,
\zeta_{42}^3y,\zeta_{42}^{18}t) \\
44 & y^2=x^3+x+t^{11} & (x,y,t)\mapsto (-x,\zeta^{11}_{44}y,
\zeta_{44}^2t) \\
66 & y^2=x^3+t(t^{11}-1) & (x,y,t)\mapsto (\zeta^2_{66}x,
\zeta_{66}^3y, \zeta_{66}^6t)\\ \hline
\end{array}
$$
\medskip

Here $\zeta_k$ denotes a primitive $k$-th root of unity.
\medskip

Now the main result of Livn\'e--Sch\"utt--Yui \cite{LSY2010} is 
to establish the modularity of these K3 surfaces.

\begin{theorem}
Let $X$ be a {\mbox{K3}} surface in the above table. Then for
each $k$, the $\ell$-adic Galois representation associated to
$T(X)$ is irreducible over $\QQ$ of dimension $\phi(k)$.
Furthermore, this $G_{\QQ}$-Galois representation is
induced from a one-dimensional Galois representation
of $\QQ(\zeta_k)$. 

All these K3 surfaces are of {\mbox{CM}} type, and are modular
(automorphic).
\end{theorem}

\proof

$\bullet$ CM type is established by realizing them
as Fermat quotients.  Since Fermat surfaces
are known to be of CM type, the result follows.

$\bullet$ Modularity (or automorphy) of the Galois
representation is established by using automorphic
induction. The restriction of the Galois representation
to the cyclotomic field $\QQ(\mu_k)$ is given
by a one-dimensional Jacobi sum Gr\"ossencharakter.
To get down to $\QQ$, we take the $\mathrm{Gal}(\QQ(\zeta_k)/\QQ)$-orbit
of the one-dimensional representation. This Galois
group has order $\phi(k)$, and we obtain the
irreducible Galois representation over $\QQ$ of dimension
$\phi(k)$.

$\bullet$ Those K3 surfaces corresponding to $k=44$ and $66$ are singular,
so their modularity has already been established by Theorem 2.2. 

\begin{remark}
When $T(X)$ is not unimodular, there are $10$ values
of $k$ such that $\rank T(X)=\phi(k)$:
$$19, 17, 13, 11, 7, 25, 5, 27, 9, 3$$
All these K3 surfaces are again dominated by
Fermat surfaces, and hence they are all of CM type.
We have also established their modularity (automorphy).
\end{remark}
\medskip

In the article of Goto--Livn\'e--Yui \cite{GLY2011},
more examples of K3 surfaces of CM type are constructed.
First we recall a classification result of Nikulin \cite{N86}.

Let $X$ be a K3 surface over $\QQ$. Let
$H^{2,0}(X)=\CC\omega_X$ where we fix a nowhere vanishing
holomorphic $2$-form $\omega_X$. Let $\sigma$ be an
involution on $X$ such that $\sigma(\omega_X)=-\omega_X$.
Let $=\Pic(X)^{\sigma}$ be the fixed part of $\Pic(X)$ by
$\sigma$.  Put $r=\mathrm{rank}\Pic(X)^{\sigma}$.
Let $T(X)_0$ be the orthogonal complement
of $\Pic(X)^{\sigma}$ in $H^2(X,\ZZ)$. Then $\sigma$ acts
by $-1$ on $T(X)_0$. 
Consider the quotient groups $(\Pic(X)^{\sigma})^*/\Pic(X)^{\sigma})$
and $(T(X)_0^*/T(X)_0)$, where $L^*$ denotes the dual lattice
of a lattice $L$.
Since $H^2(X,\ZZ)$ is unimodular,
the quotient abelian groups are canonically isomorphic
$$(\Pic(X)^{\sigma})^*/\Pic(X)^{\sigma}\simeq 
(T(X)_0^*/T(X)_0.$$
Since $\sigma$ acts as $1$ on the first quotient; and 
as $-1$ on the second quotient, this forces these 
quotient groups to be isomorphic to $(\ZZ/2\ZZ)^a$ for
some positive integer $a\in\ZZ$. 

The intersection pairing on $\Pic(X)$ induces a quadratic
form $q$ on the discriminant group with values in $\QQ$ modulo $2\ZZ$; 
we put 
$\delta=0$ if $q$ has values only in $\ZZ$, and $\delta=1$ 
otherwise.

Thus, we have a triplet of integers $(r, a, \delta)$ associated

to a K3 surface $X$ with the involution $\sigma$. 
A theorem of Nikulin \cite{N86} asserts that a pair $(X, \sigma)$ is 
classified, up to deformation, by a triplet $(r, a, \delta)$.

\begin{theorem}
There are $75$ possible triplets $(r,a,\delta)$ that 
classify pairs $(X,\sigma)$ of {\mbox{K3}} surfaces
$X$ with involution $\sigma$, up to deformation.
\end{theorem}


\begin{figure}
\centering
\includegraphics[width=10cm]{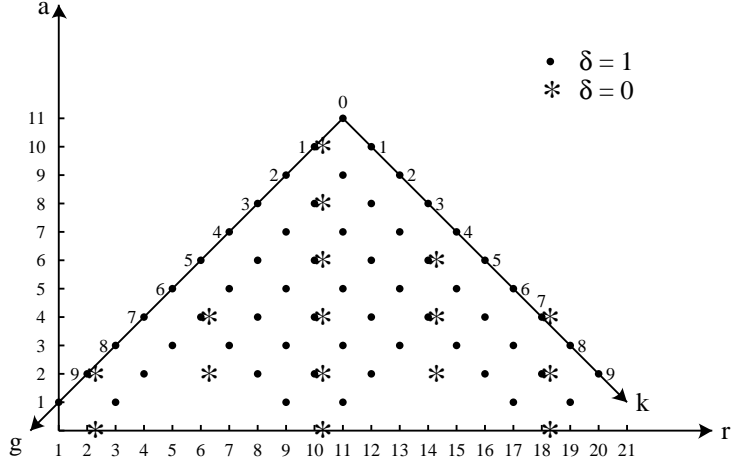}  %
\caption{Nikulin's pyramid}
\label{pyramid}
\end{figure}

Now we will realize some of these $75$ families of
K3 surfaces with involution.  We look for K3 surfaces
defined by hypersurfaces.  For this we use the famous
$95$ families of hypersurfaces in weighted
projective $3$-spaces determined by
M. Reid \cite{R79} or Yonemura \cite{Yo90}.
Let $[x_0:x_1:x_2:x_3]$ denote weighted projective
coordinates in a weighted projective $3$-space
with weight $(w_0,w_1,w_2,w_3)$.

\begin{theorem} (Goto--Livn\'e--Yui \cite{GLY2011})
Among the $95$ families of {\mbox{K3}} surfaces, all but
$9$ families of {\mbox{K3}} surfaces have an involution
$\sigma$, satisfying the following conditions:

(1) removing several monomials (if necessary) from
Yonemura's hypersurface, a new defining hypersurface 
consists of exactly four monomials, i.e., it is of Delsarte
type,

(2) the new defining hypersurface is quasi-smooth,
and the singularity configuration should
remain the same as the original defining equation
of Yonemura, 

(3) the new defining hypersurface contains
only one monomial in $x_0$ of the form $x_0^n$, 
$x_0^nx_j$, $x_i^n+x_ix_j^m$ or $x_i^nx_k+x_ix_j^m$ for some 
$j$ and $k\, (k\neq j)$ distinct from $i$. 

For $45$ (resp. $41$) {\mbox{K3}} surfaces, a defining equation 
of four monomials is of the form
$$x_0^2=f(x_1,x_2,x_3)\subset \PP^3(w_0,w_1,w_2,w_3)$$
(resp.
$$F(x_0,x_1,x_2,x_3)=0\subset \PP^3(w_0,w_1,w_2,w_3))$$
where $f$ (resp. $F$) is a homogeneous polynomial in
$x_1,x_2,x_3$ (resp. $x_0,x_1,x_2,x_3$) over $\QQ$ of degree
$\sum_{i=0}^3 w_i$.

(4) For all $86=45+41$ K3 surfaces, there is an
algorithm to compute the invariants $r$ and $a$.
\end{theorem}

\begin{remark} There are nine families for which the Theorem is not
valid. 
The three families ($\#15, \#53$ and $\#54$) do not have the required 
involution. Another different six families ($\#85,\#94, \#95$) and
($\#90, \#93, \#91$) cannot not realized as
quasi-smooth hypersurfaces in four monomials.
(We employ the numbering from Yonemura \cite{Yo90}.)
\end{remark}

\begin{proposition} (Nikulin \cite{N86}; Voisin \cite{V93})
Let $(S, \sigma)$ be a {\mbox{K3}} surface with
involution $\sigma$. Let $S^{\sigma}$ be the fixed part of $S$.
Then for $(r,a,\delta)\neq (10,10,0),(10,8,0)$,
let $S^{\sigma}=C_g\cup L_1\cup \cdots \cup L_k$ where
$C_g$ is a smooth genus $g$ curve and $L_1,\ldots, L_k$ are
rational curves. Then
$$r=11+k-g,\,\, a=11-g-k.$$

If $(r,a,\delta)=(10,10,0)$, then $S^{\sigma}=C_1\cup C_2$
where $C_i$ ($i=1,2$) are elliptic curves, and if $(r,a,\delta)=(10,8,0)$,
then $S^{\sigma}=\emptyset$.
\end{proposition}

\begin{proposition} (Goto--Livn\'e--Yui \cite{GLY2011})
Let $(S,\sigma)$ be one of the $45$ pairs in the
above theorem, which is given as the minimal
resolution of a hypersurface
$S_0\,:\,x_0^2=f(x_1,x_2,x_3)\subset\PP^3(Q)$
where $Q=(w_0,w_1,w_2,w_3)$ and $f$ is a homogeneous
polynomial of degree $\sum_{i=0}^3 w_i$. 
Let $r(Q)$ denote the number of exceptional divisors
in the resolution $S\to S_0$.
Then $r$ can be computed as follows:

(i) Suppose $w_0$ is odd. Then
$r=r(Q)-w_i+2$ if there is an odd weight $w_i\neq w_0$
such that $\gcd(w_0,w_i)=w_i\geq 2$, and $r=r(Q)+1$
otherwise.

(ii) Suppose $w_0$ is even. Then
$r=r(Q)+1-\sum_{i=1}^3 (d_i-1)(\frac{2d_i}{w_i}-1)$
where $d_i=\gcd(w_0,w_i)$.
\end{proposition}

\begin{remark}
For non-Borcea type K3 surfaces defined by equations
of the form $x_0^2x_i=f(x_1,x_2,x_3)$ for some $i\in\{1,2,3\}$,
the invariants $r$ and $a$ (or equivalently, $g$ and $k$) are
also computed.
\end{remark}

\begin{proposition} (Goto--Livn\'e--Yui \cite{GLY2011})
At least the $39$ triplets $(r,a,\delta)$ of integers are realized
by the $86$ families of {\mbox{K3}} surfaces.
\end{proposition}

\begin{theorem} (Goto--Livn\'e--Yui \cite{GLY2011})
For the $86$ families of {\mbox{K3}} surfaces, there are 
subfamilies of {\mbox{K3}} surfaces which are of {\mbox{CM}} type, 
that is, there is a {\mbox{CM}} point in the moduli space of 
{\mbox{K3}} surfaces.  Indeed, at each {\mbox{CM}} point, 
the {\mbox{K3}} surface is realized as a quotient of a Fermat 
(or Delsarte) surface. Consequently, it is modular (automorphic) 
by automorphic induction.
\end{theorem}

\begin{remark}
Nikulin's Theorem 2.9 gives in total $75$ triplets of integers $(r, a,\delta)$
that classify the isomorphism classes of pairs $(S, \sigma)$ of K3 surfaces
$S$ with non-symplectic involution $\sigma$. With our choice of K3 surfaces
of CM type, we can realize at least $39$ triplets.  This is based on our
calculations of the invariants $r$ and $a$. Since we have not yet computed
the invariant $\delta$ for the $39$ cases, this number may increase somewhat.
\end{remark}

\begin{remark}
It is known (Borcea \cite{B97}) that over $\CC$ the moduli spaces
of Nikulin's K3 surfaces are arithmetic quotients of type IV (Shimura
varieties).
Recently, Ma \cite{Ma} has shown the rationality of the moduli
spaces for the $67$ triplets.  Our result above
gives one explicit CM point in these moduli spaces. 
CM type implies that these families must be isolated points in the
moduli space. We are not able to show the denseness of CM points, however.
\end{remark}

It is notoriously difficult to compute Hodge groups,
and the above theorem implies the commutativity
of the Hodge group.

\begin{corollary}
The Hodge groups of these $86$ {\mbox{K3}} surfaces are
all commutative, i.e., copies of $\GG_m$'s over $\CC$.
\end{corollary}

\begin{proposition} (Goto-Livn\'e--Yui \cite{GLY2011})
We classify the $86$ hypersurfaces into four types:

\begin{itemize}
\item $45$ families have the form:
$x_0^2=f(x_1,x_2,x_3)$ and $\sigma(x_0)=-x_0$. \\
\item $33$ families cannot be put in the
form (1) but defined by a hypersurface of the form 
$F(x_0,x_1,x_2,x_3)=0$ with $4$ monomials, and $\sigma$
can be described explicitly.
\item $4$ families, 
after changing term the $x_0^3$ to $x_0^2x_1$ and then removing
several terms, can be put into the form $F(x_0,x_1,x_2,x_3)=0$ of 
four monimials equipped with an explicit involution $\sigma$.
\item The remaining $4$ families can be put in
the form $F(x_0,x_1,x_2,x_3)=0$ of four monomials equipped with a different
kind of involution.  
\end{itemize}
\end{proposition}

\begin{example}

(1) Here is one of the $45$ cases,
$\#78$ in Yonemura's list. The weight is
$(11,6,4,1)$, and Yonemura's hypersurface is
$$x_0^2=x_1^2x_2+x_1^3x_3^4+x_1x_2^4+x_2^5x_3^2+x_3^{22}.$$
We can remove
$x_1^3x_2^4$ and $x_2^5x_3^2$ so the new hypersurface is
$x_0^2=x_1^2x_2+x_1x_2^4+x_3^{22}$. The singularity is
of type $A_1+A_3+A_5$.

(2) Here is one of the second cases, 
$\#19$ in Yonemura's list. The weight is
$(3,2,2,1)$ and Yonemura's hypersurface is
$$F(x_0,x_1,x_2,x_3)=x_0^2x_1+x_0^2x_2+x_0^2x_3^2+x_1^4+x_2^4+x_3^4.$$
We can remove $x_0^2x_1$ or $x_0^2x_2$, $x_0^2x_3^2$,
so the new hypersurface is
$x_0^2x_1+x_1^4+x_2^4+x_3^4.$
The singularity is of type $4A_1+A_2$.

(3) Here is one of the third cases, $\#18$ in Yonemura's list. This 
hypersurface $x_0^3+x_1^3+x_0x_2^3+x_1x_2^3+x_2^4x_3+x_3^9$ acquires 
an involution if we
replace $x_0^3$ by $x_0^2x_1$ and remove the terms $x_0x_2^3$ and $x_1x_2^3$.

(4) The last case is $\#52$ in Yonemura's list.

(5) $\#95$ in Yonemura's list. This hypersurface
has an involution, but cannot be realized by
a quasi-smooth hypersurface with four monomials.
\end{example}
\medskip

\subsection{The modularity of higher dimensional
Galois representations arising from Calabi--Yau
threefolds over $\QQ$}
\medskip

There are several new examples of modular 
non-rigid Calabi--Yau threefolds $X$ over $\QQ$ with $B_3\geq 4$.
These examples were constructed after the article Yui\cite{Yui2003},
some of which have already been discussed in the article
of E. Lee \cite{L2008}.

There are several approaches (with non-empty intersection)
to produce these new modular examples: 

{\bf (1)} Those non-rigid Calabi--Yau threefolds $X$ over $\QQ$ such that
the semi-simplification of $H^3_{et}(\ov{X},\QQl)$ is highly reducible
and splits into smaller dimensional irreducible Galois representations. 
For instance, most known cases are when the third cohomology
group splits into $2$-dimensional or $4$-dimensional pieces. (Examples
of E. Lee, Hulek and Verrill, Sch\"utt, Cynk and Meyer, and a
recent example of Bini and van Geemen, and Sch\"utt, and more.)
The article of E. Lee \cite{L2008} gives reviews on the modularity
of non-rigid Calabi--Yau threefolds over $\QQ$, up to 2008.

{\bf (2)} Those Calabi--Yau threefolds $X$ over $\QQ$
such that the the $\ell$-adic Galois representations arising
from $H^3(\ov{X},\QQl)$ are irreducible and have small dimensions (e.g.,
$4, 6,$ or $8$.)
(Examples of Livn\'e--Yui, Dieulefeit--Pacetti--Sch\"utt, and more.)

{\bf (3)} Those Calabi--Yau threefolds over $\QQ$ which
are of CM type, thus the Galois representations
are induced by one-dimensional representations. 
(Examples of Rohde, Garbagnati--van Geemen, Goto--Livn\'e--Yui, and more.)

{\bf (4)} Given a Calabi--Yau threefold $Y$ over $\QQ$, if we
can construct an algebraic correspondence defined over $\QQ$
to some modular Calabi--Yau threefold $X$ over $\QQ$, the Tate
conjecture asserts that their $L$-series should coincide.
This will establish the modularity of $Y$. 
\medskip

\noindent{\bf (1): Calabi--Yau threefolds in the category (1), i.e, with
highly reducible Galois representations} 

The general strategy is to consider Calabi--Yau
threefolds which contain elliptic ruled
surfaces.  This is formulated by Hulek and Verrill
\cite{HV2005a}.

\begin{proposition} (Hulek and Verrill \cite{HV2005a})
Let $X$ be a Calabi--Yau threefold over $\QQ$. Suppose
that $X$ contains birational ruled elliptic surfaces 
$S_j,\,j=1,\dots, b$ over $\QQ$ and whose cohomology classes span 
$H^{2,1}(X)\oplus H^{1,2}(X)$ (so $b=h^{2,1}(X)$.)
Let $\rho$ be the $2$-dimensional Galois representation
given by the kernel $U$ from the exact sequence
$$0\to U\to H^3_{et}(\ov{X},\QQl)\to \oplus H^3_{et}(\ov{S}_j,\QQl)\to 0.$$
Then $X$ is modular, that is,
$$L(X,s)=L(f_4,s)\prod_{j=1}^b L(g_2^j,s-1)$$
where $f_4$ is a weight $4$ modular form associated to
$\rho$ and $g_2^j$ are the weight $2$ modular forms
associated to the base elliptic curves $E_j$ of the birational
ruled surfaces $S_j$.
\end{proposition}

The requirement that the third cohomology group splits
as in the proposition is rather restrictive. Several
examples of Calabi--Yau threefolds satisfying this
condition are given by Hulek and Verrill \cite{HV2005a},
E. Lee \cite{L2006}, 
and Sch\"utt \cite{S2006}. 
They are constructed as
resolutions of fiber products of semi-stable rational 
elliptic surfaces with section.
Another series of examples along this line 
due to Hulek and Verrill \cite{HV2005b}
are toric Calabi--Yau threefolds associated to the
root lattice $A_4$. 

Cynk and Meyer \cite{CM2008} have established the modularity of
$17$ nonrigid double octic Calabi--Yau threefolds over $\QQ$ with
$B_3=6$. The Galois representations decompose into $2$- and $4$-dimensional
sub-representations such that the $L$-series of each such sub-representation
is of the form $L(g_4,s), \, L(g_2,s-1)$ or $L(g_2\times g_3,s)$, where
$g_k$ is a weight $k$ cusp form.
\medskip

Now we will consider another construction due to E. Lee \cite{L2006, L2011}  
of Calabi--Yau threefolds associated to the Horrocks--Mumford
vector bundle of rank $2$.  It is well known that
Horrocks--Mumford quintics are determinantal
quintics. The Schoen quintic $Q: \sum_{i=0}^4x_i^5-5\prod_{i=0}^4x_i=0$
is the early example of this type. 

Lee has constructed more Calabi--Yau threefolds.
Let $y\in\PP^4$ be a generic point, and
define the matrices
$$M_y(x)=\left(\begin{array}{ccccc}
      x_0y_0 & x_3y_2 & x_1y_4 & x_4y_1 & x_2y_3 \\
     x_3y_3 & x_1y_0 & x_4y_2 & x_2y_4 & x_0y_1 \\
     x_1y_1 & x_4y_3 & x_2y_0 & x_0y_2 & x_3y_4 \\
     x_4y_4 & x_2y_1 & x_0y_3 & x_3y_0 & x_1y_2 \\
     x_2y_2 & x_0y_4 & x_3y_1 & x_1y_3 & x_4y_0\end{array}\right)$$
and
$$L_y(z)=\left(\begin{array}{ccccc}
     z_0y_0 & z_2y_4 & z_4y_3 & z_1y_2 & z_3y_1 \\
     z_4y_1 & z_1y_0 & z_3y_4 & z_0y_3 & z_2y_2 \\
     z_3y_2 & z_0y_1 & z_2y_0 & z_4y_4 & z_1y_3 \\
     z_2y_3 & z_4y_2 & z_1y_1 & z_3y_0 & z_0y_4 \\
     z_1y_4 & z_3y_3 & z_0y_2 & z_2y_1 & z_4y_0\end{array}\right)$$
Note that $M_y(x)z=L_y(z)x$. Then
$$X_y:=\{\,\mathrm{det}\,M_y(x)=0\}\subset\PP^4(x)$$
and
$$X^{\prime}_y:=\{\,\mathrm{det}\, L_y(z)=0\}\subset\PP^4(z)$$
are Horrocks--Mumford quintics.  Define a threefold
$\tilde X_y$ in $\PP^4(x)\times\PP^4(z)$ as a common
partial resolution of $M_y(x)z=0.$
We ought to know the singularities of $\tilde X_y$.
Let $\PP^2_+:=\{\,y: y_1-y_4=y_2-y_3=0\,\}$. Given
$y\in\PP^2_+$, the point $(1:0:0:0:0)$ is a
singular point over $\CC$ of $X_y$ if and only if one of the
coordinates of $y$ is zero. For $y=(0:1:-1:-1:1:1),
(0:2:3\pm\sqrt{5}:3\pm\sqrt{5}:2),
(2:-1:0:0:-1),(2:\pm\sqrt{5}-1:0:0:\pm\sqrt{5}-1)\in\PP^2_+$,
$X_y$ contains the Heisenberg orbits of $(1:0:0:0:0)$ and
$(1:1:1:1:1)$ as nodes over $\CC$.

\begin{proposition} (Lee \cite{L2006, L2008, L2011})
(a) Let $y=(0:1:-1:-1:1)$, and write
$X_{(0:1:-1:-1:1)}=X$ for short. Then the Calabi--Yau threefold
$\hat X$ obtained by crepant resolution of singularities has
$B_3=6$. Furthermore,  
$$H^3_{et}(\ov{\hat X}, \QQl)=V\oplus
H^2(\ov{S},\QQl)(-1)$$
where $V$ is $2$-dimensional and associated to the modular form of weight
$4$ and level $55$, and $H^2(\ov{S},\QQl)$ is $4$-dimensional
and is isomorphic to 
${\mbox{Ind}}^{G_{\QQ}}_{G_{\QQ(i)}} H^1(\ov{E},\QQl)(-1)$
where $E$ is an elliptic curve over $\QQ(i)$ coming from $E$ over $\QQ$. 

(b) Let $y=(2:-1:0:0:-1)$, and let $X_{(-2:-1:0:0:-1)}=X^{\prime}$
for short. Then the Calabi--Yau threefold $\hat{X}^{\prime}$ obtained
by crepant resolution of singularities has $B_3=4$.
Furthermore  
$$H^3_{et}(\ov{\hat X},\QQl)
=V\oplus H^1(\ov{E}_2, \QQl)(-1)$$ where
$V$ is the $2$-dimensional Galois representation 
associated to the modular form of weight $4$ of level $55$, 
and $H^1(\ov{E}_2,\QQl)$ is associated with the modular form
of weight $2$ and level $550$. 
Furthermore, the $L$-series of $\hat X$ is given, up to Euler
factors at the primes of bad reduction, by
$$L(\hat X,s)=L(f,s)L(g, s-1)$$
where $f$ is the unique normalized cusp form of weight $4$
and level $55$ and $g$ is the normalized cusp form of weight $2$
and level $550$.

(c) Now consider a smooth (big) resolution $X$ of the 
$(\ZZ/2\ZZ)$-quotient of the Schoen quintic 
$Q: x_0^5+x_1^5+x_2^5+x_3^5+x_4^5-5x_0x_1x_2x_3x_4=0$.
(There is the $(\ZZ/2\ZZ)$-action on $Q$ induced by the involution
on $\PP^4$ defined by $\iota[x_0:x_1:x_2:x_3:x_4]=[x_0:x_4:x_3:x_2:x_1]$.)
Then $X$ is a Calabi--Yau threefold defined over $\QQ$
with $B_3=4$. The Calabi--Yau threefold $X$ is modular,
and up to Euler factors at primes of bad reduction at $p=2$ and $5$,
$L(X,s)$ is given by
$$L(X,s)=L(f,s)L(g,s-1)$$
where $f$ is the unique normalized cusp form of weight $4$ and level
$25$, and $g$ is a weight $2$ cusp form of level $50$.
\end{proposition}

A recent example due to Bini and van Geemen \cite{BvG2011},
and Sch\"utt \cite{S2011} is the Calabi--Yau threefold
called Maschke's double octic, which arises as the double covering 
of ${\PP}^3$ branched along Maschke's surface $S$.
The Maschke octic surface $S$ is
defined by the homogeneous equation 
$$S=\large{\sum_{i=0}^3 x_i^8+14\sum_{i<j} x_i^4x_j^4
+168 x_0^2x_1^2x_2^2x_3^2=0\large}\subset\PP^3$$ 

Now let $X$ be the double cover of $\PP^3$ along $S$. This is
a smooth Calabi--Yau threefold defined over $\QQ$.
Let $Y$ be the desingularization of the quotient of $X$ by a 
suitable Hisenberg group. Then $Y$ is also a smooth Calabi--Yau
threefold defined over $\QQ$. The results of Bini and van Geemen,
and Sch\"utt, are summarized in the following

\begin{proposition} {\sl (a) The Maschke double octic
Calabi--Yau threefolds $X$ is modular over $\QQ$. The
third cohomology group $H^3_{et}(\ov{X},\QQl)$ has $B_3(X)=300$. The Galois
representation of $H^3_{et}(\ov{X},\QQl)$ decomposes completely over $\QQ(i)$ 
into $2$--dimensional Galois representations which descend to
$\QQ$, and the latter correspond to modular forms of weight $4$, or 
modular forms of weight $2$.

(b) Let $Y$ be the desingularization of $X$ by a suitable
Heisenberg group. Then $Y$ is modular over $\QQ$. The
third cohomology group $H^3_{et}(\ov{Y},\QQl)$ has $B_3(Y)=30$. The Galois
representation of $H^3_{et}(\ov{Y},\QQl)$ decomposes completely over $\QQ$ 
into $2$-dimensional Galois representations, and the latter correspond 
to modular forms of weight $4$ and modular forms of weight $2$.

(c) The Maschke surface $S$ has Picard number $\rho(S)=202$.
The second cohomology group $H^2(\ov{S},\QQl)$ has $B_2(S)=302$. The
Galois representation of the transcendental part has dimension $100$,
which splits into $2$ or $3$-dimensional Galois sub-representations
over $\QQ$, and the latter correspond to modular forms
of weight $3$, or modular forms of weight $2$.}
\end{proposition}
\medskip

\noindent{\bf (2) Calabi--Yau threefolds in the category (2), i.e,
with irreducible Galois representations}

Consani--Scholten \cite{CS01} constructed a Calabi--Yau
threefold over $\QQ$ as follows. Consider the Chebyshev polynomial
$$P(y,z)=(y^5+z^5)-5yz(y^2+z^2)+5yz(y+z)+5(y^2+z^2)-5(y+z)$$
and define an affine variety $X$ in $\mathbf{A}^4$ by
$$X: P(x_1,x_2)=P(x_3,x_4)$$
and let $\overline{X}\subset \PP^4$ be its projective
closure. Then $\overline{X}$ has $120$ ordinary double points.
Let $\tilde{X}$ be its small resolution. Then $\tilde{X}$ 
is a Calabi--Yau threefold with $h^{1,1}(\tilde{X})=141$ and $h^{2,1}(\tilde{X})=1$.
$\tilde{X}$ is defined over $\QQ$ and the primes $2,3,5$ are bad primes.
$H^3(\ov{\tilde{X}},\QQl)$ gives rise to a $4$-dimensional $\ell$-adic Galois
representation, $\rho$, which is irreducible
over $\QQ$. Let $F=\QQ(\sqrt{5})$, and let
$\lambda\in F$ be a prime above $\ell$. Then the restriction 
$\rho|_{\mbox{Gal}(\bar{\QQ}/F)}$ is reducible as a representation
to $GL(4,F_{\lambda})$: There is a Galois representation
$\sigma : \mbox{Gal}(\bar{\QQ}/F)\to GL(2, F_{\lambda})$ such that
$\rho=\mbox{Ind}^{\QQ}_F\sigma$. Consani--Scholten conjectured
the modularity of $\tilde{X}$.

\begin{theorem} (Dieulefait--Pacetti--Sch\"utt \cite{DPS2010})
The Consani--Scholten Calabi--Yau threefold $\tilde{X}$ over
$\QQ$ is Hilbert modular. That is, the $L$-series associated
to $\sigma$ coincides with the $L$-series of a Hilbert modular
newform $\mathfrak{f}$ on $F$ of weight $(2,4)$ and conductor
$\mathfrak{c}_{\mathfrak{f}}=(30)$.
\end{theorem}
\medskip

The first example of Siegel modular varieties, as moduli spaces, of
Calabi--Yau threefolds 
was given by van Geemen and Nygaard \cite{GN95}.
Recently, a series of articles by Freitag and Salvati Manni \cite{FM} 
on Siegel modular threefolds
which admit Calabi--Yau models have appeared. 
The starting point is the van Geemen--Nygaard rigid Calabi--Yau
threefold defined by a complete intersection $Y$
of degree $(2,2,2,2)$ in $\PP^7$ by the equations:
$$
\begin{array}{ccccccccc}
Y_0^2&=& X_0^2&+& X_1^2&+& X_2^2&+& X_3^2 \cr
Y_1^2&=& X_0^2&-& X_1^2&+& X_2^2&-& X_3^2 \cr
Y_2^2&=& X_0^2&+& X_1^2&-& X_2^2&-& X_3^2 \cr
Y_3^2&=& X_0^2&-&X_1^2&-& X_2^2&+& X_3^2\cr
\end{array}
$$
A smooth small resolution of $Y$, denoted by $X$, is a
Calabi--Yau threefold with $h^{1,1}(X)=32$
and $h^{2,1}(X)=0$, so $X$ is rigid and hence
is modular, indeed, the $L$-series of the Galois
representation associated to $H^3(X,\QQl)$ is determined
by the unique weight $4$ modular form on $\Gamma_0(8)$. 

Several examples of non-rigid Calabi--Yau threefolds
$Z$ are obtained by quotienting $X$ by finite group
actions. Many of the resulting varieties $Z$ are Calabi--Yau
threefolds with small third Betti number. 

\begin{remark}
For instance, Freitag--Salvati Manni has constructed such
Calabi--Yau threefolds $Z$. The Galois representations associated to 
$H^3(Z,\QQl)$ of their Calabi--Yau threefolds $Z$ should
be studied in detail.  They should decompose into direct sum of
those coming from the rigid Calabi--Yau threefold $X$, and 
those coming from elliptic surfaces or surfaces of 
higher genus arising from fixed points of finite groups in question.  

In particular, this would imply that proper
Siegel modular forms will not arise from
these examples. So far as I know, we do not have 
examples of Calabi--Yau threefolds over $\QQ$ with $B_3(X)=4$
whose $L(X,s)$ comes from a Siegel modular form on $Sp_4(\ZZ)$
or its subgroups of finite index.
\end{remark}
\medskip

\noindent{\bf (3) Calabi--Yau threefolds of CM type}

We next consider Calabi--Yau threefolds which
we will show to be of CM type.  Then
the Galois sub-representations associated to the
third cohomology groups 
are induced by one-dimensional ones. Then, by 
applying the automorphic induction process, we will
establish the automorphy of Calabi--Yau threefolds.
These Calabi--Yau threefolds are realized
as quotients of products of K3 surfaces and
elliptic curves by some automorphisms. 
Rohde \cite{R2008}, Garbagnati--van Geemen \cite{GG2010},
and Goto--Livn\'e--Yui \cite{GLY2011} produce examples of CM type 
Calabi--Yau threefolds with this approach.
\medskip

Let $S$ be a K3 surface with an involution $\sigma$
acting on $H^{0,2}(S)$ by $-1$ discussed in
subsection 2.2. Let $E$ be an elliptic curve with
the standard involution $\iota$. Consider the quotient
of the product $E\times S/\iota\times\sigma$. This
is a singular Calabi--Yau threefold having only cyclic
quotient singularities. Resolving singularities
we obtain a smooth crepant resolution $X$. Since
the invariants of $X$ are determined by a triplet of integers $(r,a,\delta)$
associated to $S$, we will write $X$ as $X(r,a,\delta)$. 
The Hodge numbers and the Euler characteristic 
of $X$ depend only on $r$ and $a$:
$$h^{1,1}(X)=5+3r-2a,\, h^{2,1}(X)=65-3r-2a$$
and
$$ e(X)=2(h^{1,1}(X)-h^{2,1}(X))=6(r-10).$$

\begin{theorem} (Goto--Livn\'e--Yui \cite{GLY2011})
Let $(S, \sigma)$ be one of the {\mbox{K3}} surfaces
defined over $\QQ$ in Theorem 2.9. Let $E$ be an 
elliptic curve over $\QQ$ with the standard involution 
$\iota$.  Let $X=X(r,a,\delta)$ be a smooth Calabi--Yau
threefold. Then the following assertions hold:

(1) $X$ is of {\mbox{CM}} type if and only if $E$ is of
{\mbox{CM}} type.

(2) If $X$ is of {\mbox{CM}} type, then the Jacobian variety
$J(C_g)$ of $C_g$ in $S^{\sigma}$ is
also of {\mbox{CM}} type, provided that the K3 surface component
is of the form $x_0^2=f(x_1,x_2,x_3)$ with involution $\sigma(x_0)=-x_0$.

(3) $X$ is modular (automorphic).
\end{theorem}

Sketch of Proof: (1) The Hodge structure $h_X$ of type
$(3,0)$ of $X$ is given by the tensor product $h_S\otimes h_E$
 of the Hodge structures $h_S$ of type $(2,0)$ and $h_E$ of
type $(1,0)$.  
Then $h_X$ is of CM type if and only if both
$h_S$ and $h_E$ are of CM type. Since $S$ is already
of CM type, we only need to require that $E$ is
of CM type.

(2) When the K3 surface is defined by a hypersurface
of the form $x_0^2=f(x_1,x_2,x_3)$ and the involution
$\sigma$ takes $x_0$ to $-x_0$, then the curve $C_g$ in 
the fixed locus $S^{\sigma}$ is obtained by
putting $x_0=0$, and hence it is also 
of Delsarte type. Hence the Jacobian variety $J(C_g)$ of
$C_g$ is also of CM type.

{\it When the hypersurface defining the K3 surface is not
of the above form and the involution is more complicated, we
ought to check each case whether $C_g$ is of Delsarte
type or not.}

(3) $S$ is modular by Theorem 2.11, and $E$ is modular
by Wiles et al. Hence $X$ is modular (automorphic).
\medskip

Now we will discuss mirror Calabi--Yau threefolds
of $X=X(r,a,\delta)$. 

\begin{proposition} (Borcea \cite{B97} and Voisin \cite{V93})
Given a Calabi--Yau threefold $X=X(r,a,\delta)$, there
is a mirror Calabi--Yau threefold $X^{\vee}$
such that $X^{\vee}$ is realized as a crepant
resolution of a quotient of $E\times S/\iota\times\sigma$
and $X^{\vee}$ is characterized by the invariants
$(20-r, a, \delta)$. The Hodge numbers of $X^{\vee}$ are
$$h^{1,1}(X^{\vee})=5+3(20-r)-2a=65-3r-2a=h^{2,1}(X),$$
$$h^{2,1}(X^{\vee})=65-3(20-r)-2a=5+3r-2a=h^{1,1}(X)$$
and the Euler characteristic is
$$e(X^{\vee})=-12(r-10)=-e(X).$$

In terms of $g$ and $k$, $r=11-g+k,\,a=11-g-k$, and 
$$h^{1,1}(X)=1+r+4(k+1),\,\, h^{2,1}(X)=1+(20-r)+4g$$
and the Euler characteristic is
$e(X)=12(1+k-g).$
\end{proposition}

\begin{remark}
Mirror symmetry of Calabi--Yau threefolds of $\mbox{K3}\times E$
do come from mirror symmetry of K3 surfaces. 
Mirror symmetry for K3 surfaces is the correspondence
$r \leftrightarrow (20-r)$. In fact, one can see this in Figure 1: Nikulin's
pyramid.  Given a K3 surface $S$,
we try to look for a mirror K3 surface $S^{\vee}$ satisfying
this correspondence.  This correspondence is established at
a special CM point in the moduli space.

We know that the $95$ families of K3 surfaces of
Reid and Yonemura are not closed under mirror
symmetry. Only $54$ families of K3 surfaces
with involution have mirror partners within the $95$ families.

A recent article of
Artebani--Boissi\`ere--Sarti \cite{ABS2011} also considers
this type of mirror symmetry for K3 surfaces.
\end{remark}

\begin{proposition}
For a mirror $X^{\vee}$, the assertion is true
if the mirror {\mbox{K3}} surface is also of Delsarte type.
\end{proposition}

\begin{proposition} 
Let $(S,\sigma)$ be one of the $86$ families of
{\mbox{K3}} surfaces with involution $\sigma$.  Let $X=X(r,a,\delta)$
and $X^{\vee}=X(20-r,a,\delta)$ be mirror
pairs of Calabi--Yau threefolds defined
over $\QQ$, where we suppose that the {\mbox{K3}} surface component
is of Delsarte type.  Then $X$ and $X^{\vee}$
have the same properties: 
\begin{itemize}
\item  $X$ is of {\mbox{CM}} type
if and only if $X^{\vee}$ is of {\mbox{CM}} type, and
\item $X$ is modular if and only if $X^{\vee}$ is modular.
\end{itemize}
\end{proposition}

Rohde \cite{R2008}, Garbagnati and van Geemen \cite{GG2010} 
and Garbagnati \cite{G2010} constructed Calabi--Yau 
threefolds which are quotients of the products of 
K3 surfaces and
elliptic curves by non-symplectic automorphisms 
of higher order (than $2$), that is, order $3$, or
order $4$. These Calabi--Yau threefolds are parametrized
by Shimura varieties.

Rohde \cite{R2008} constructed families of Calabi--Yau
threefolds as the desingularization of the quotient
$S\times E$ by an automorphism of order $3$ where $E$
is the unique elliptic curve with an automorphism
$\alpha_E$ of order $3$, and $S$ is a K3
surface with an automorphism $\alpha_S$ of order $3$ which
fixes $k$ rational curves and $k+3$ isolated points for some
integer $k, \, 0\leq k\leq 6$.
\medskip

Let $\xi$ be a primitive
cube root of unity.  Choose the specific elliptic curve: 
$E=\CC/\ZZ+\xi\ZZ$, which has a Weierstrass model 
$y^2=x^3-1$ and $\alpha_E: (x,y)\mapsto (\xi x, y)$. 
Also choose some specific K3 surface $S$. Let
$$S=S_f: Y^2=X^3+f(t)^2,\, f=gh^2,\, \mbox{deg}(f)=6$$
where $t$ is the coordinate on $\PP^1$, and $f$ has four
distinct zeros. 
$S_f$ has an automorphism of order $3$:
$\alpha_f: (X,Y,t)\mapsto (\xi X, Y, t)$. 
Now define a Calabi--Yau threefold $X_f$ as the desingularization
of $S_f\times E$ by the automorphism $\alpha:=\alpha_f\times\alpha_E $.
Note that $S_f$ is birationally isomorphic to $(C_f\times E)/(\beta_f\times
\alpha_E)$ where $C_f: v^3=f(t)$ and $\beta_f: C_f\to C_f,\, 
(t,v)\mapsto (t,\xi v)$.  Then $X_f$ is birationally isomorphic to
$$(C_f\times E\times E)/H\quad\mbox{where}\quad
H=<\beta_f\times \alpha_E\times 1_E, 1_{C_f}\times \alpha_E\times
\alpha^{-1}_E>.$$

Rohde \cite{R2008} worked out the case $\mbox{deg}(g)=\mbox{deg}(h)=2$.
Garbagnati--van Geemen \cite{GG2010} considered the other cases,
i.e., $\mbox{deg}(g)=4,\mbox{deg}(h)=1$ and $\mbox{deg}(g)=6,\mbox{deg}(h)=0$.
The Hodge numbers of the Calabi--Yau threefold $X_f$ are computed and
the results are tabulated as follows:
\medskip

$$\begin{array}{|c|c|c|c|c|c|} \hline
\mbox{deg}(g) & \mbox{deg}(h) & g(C_f) & h^{2,1}(X_f) & h^{1,1}(X_f)& k \\ \hline
6 & 0 & 4 & 3 & 51 & 3 \\ \hline
4 & 1 & 3 & 2 & 62 & 4 \\ \hline
2 & 2 & 2 & 1 & 73 & 5 \\ \hline
0 & 3 & 1 & 0 & 84 & 6 \\ \hline
\end{array}
$$
\medskip

For $k>2$, the Calabi--Yau threefold is constructed by considering
the curve $C_{\ell} : v^6=\ell(t),\,\mbox{deg}(\ell)=12$ such
that $\ell(t)$ has five double zeros. It has the order $3$
automorphism $\beta_{\ell} : (t,v)\mapsto (t,\xi v)$. The quotient
$(C_{\ell}\times E)/(\beta_{\ell}\times \alpha_E)$ is a K3 surface
$S_{\ell}$ which has the elliptic fibration with Weierstrass 
equation $Y^2=X^3+\ell(t)$. Then the desingularization of
$(S_{\ell}\times E)/(\alpha_{\ell}\times\alpha_E)$ is a Calabi--Yau
threefold $X_{\ell}$ with $h^{2,1}(X_{\ell})=4$ and 
$h^{1,1}(X_{\ell})=40$.

Similarly, for $k=1$, one can find a Calabi--Yau threefold $X$ with
$h^{2,1}(X)=5$ and $h^{1,1}(X)=29$.

For details of these two cases, see Garbagnati--van Geemen
\cite{GG2010} or Rohde \cite{R2008}.
\medskip

We summarize the above discussion in the following form.

\begin{proposition} The Calabi--Yau threefold $X_f$ (resp. $X_{\ell}$) 
constructed above is of {\mbox{CM}} type if and only if the
Jacobian variety $J(C_f)$ (resp.
$J(C_{\ell})$) is of {\mbox{CM}} type. In this case, $X_f$ (resp. $X_{\ell})$
is modular (automorphic).  
\end{proposition}
\medskip

We will also mention results of Garbagnati
\cite{G2010}, which generalize the method of 
Garbagnati--van Geemen \cite{GG2010} to automorphisms of order
$4$. In order to construct Calabi--Yau
threefolds with a non-symplectic automorphism
of order $4$, start with the hyperelliptic
curves 
$$C_{f_g}\,:\, z^2=tf_g(t^2),\,\deg(f_g)=g,\quad 
\mbox{$f_g$ without multiple roots}.$$
$C_{f_g}$ has the automorphism $\alpha_C: (t,z)\mapsto 
(-r,iz)$. We consider the cases $g=2$ or $3$.
Let $E_i : v^2=u(u^2+1)$ be the elliptic curve and
let $\alpha_E: (u,v)\mapsto (-u, iv)$ be the autmorphism
of $E_i$. 
Now take the quotient of the product
$E_i\times C_{f_g}/\alpha_E\times \alpha_C$. Then
the singularities are of A-D-E type, and the desingularization
defines a K3 surface, $S_{f_g}$, with the automorphism
$\alpha_S: (x,y,s)\mapsto (-x,iy,s)$. In fact, there is
the elliptic fibration
$\mathcal{E}: y^2=x^3+xsf_g(s)^2$ and the map
$\pi : E_i\times C_{f_g}\to\mathcal{E}$ defined by
$((u,v);(z,t))\mapsto (x:=uz^2, y:=vz^3, s:=t^2)$ is
the quotient map 
$E_i\times C_{f_g}\to (E_i\times C_{f_g})/\alpha_E\times\alpha_C$. 
The K3 surface $S_{f_g}$ thus obtained has large Picard number. 
Indeed, $\rank(T_{S_{f_g}})\leq 4$ if $g=2$ and $\leq 6$ if $g=3$.  
Also, $S_{f_g}$ admits the order $4$ non-symplectic automorphism $\alpha_S$, 
and the fixed loci of $\alpha_S$ and of $\alpha_S^2$ contain
no curves of genus $>0$.  The K3 surface $(S_{f_g},\alpha_S^2)$ with 
involution $\alpha_S^2$ indeed corresponds to the triplet $(18,4,1)$ for $g=2$
and $(16,6,1)$ for $g=3$.

Once we have a family of K3 surfaces with non-symplectic
automorphism $\alpha_S$ of order $4$, we can construct
a family of Calabi--Yau threefolds as the quotient of
the product of $S_{f_g}$ with the elliptic curve $E_i$.  

\begin{proposition} (Garbagnati \cite{G2010})
There is a desingularization $Y_{f_g}$ of

\noindent $(E_i\times S_{f_g})/(\alpha_E^3\times \alpha_S)$
which is a smooth Calabi--Yau threefold with
$$h^{1,1}(Y_{f_g})=\begin{cases} 73 & \mbox{if $g=2$} \\
                                56 & \mbox{if $g=3$}\end{cases}\quad\mbox{and}\quad
h^{2,1}(Y_{f_g})=\begin{cases} 1 & \mbox{if $g=2$} \\
                               2 & \mbox{if $g=3$}\end{cases}.$$  
\end{proposition}

\begin{proposition}
The Calabi--Yau threefold $Y_{f_g}$ is of {\mbox{CM}} type if
and only if the Jacobian variety $J(C_{f_g})$ is
of {\mbox{CM}} type.  When it is of {\mbox{CM}} type, $Y_{f_g}$ is automorphic.
\end{proposition}

\section{The modularity of mirror maps of
Calabi--Yau varieties, and mirror moonshine} 

$\bullet$ {\bf Modularity of solutions of
Picard--Fuchs differential equations}

Let $n\in\NN$ and let
$M_n:=U_2\oplus (-E_8)^2\oplus <-2n>$ be a lattice
of rank $19$. Here $U_2$ is the usual hyperbolic lattice of
rank $2$ and $-E_8$ is the unique negative-definite
unimodular lattice of rank $8$, and $<-2n>$ denotes the
rank $1$ lattice $\ZZ v$ with its bilinear form determined by $<v,v>=-2n$.  
We consider a one-parameter family of $M_n$ polarized K3 surfaces
$X_t$ over $\QQ$ with generic Picard number $\rho(X_t)=19$.
Here by a $M_n$-polarized K3 surfaces, we mean K3 surfaces $X_t$
such that $T(X_t)$ is primitively embedded into $U_2\oplus\ZZ u$ 
where $u$ is a vector of height $2n,\, n\in\NN$.
The Picard--Fuchs differential equation of $X_t$ is of order $3$.
It is shown by Doran \cite{D20} that for such a family of
K3 surfaces $X_t$, there is a family of elliptic curves $E_t$ such
that the order $3$ Picard--Fuchs differential equation of $X_t$
is the symmetric square of the order $2$ differential equation
associated to the family of elliptic curves $E_t$. The existence
of such a relation stems from the so-called Shioda--Inose
structures of $X_t$ (or by Dolgachev's result \cite{Dol} which
asserts that the coarse moduli space of $M_n$-polarized
K3 surfaces is isomorphic to the moduli space of
elliptic curves with level $n$ structure). 
Long \cite{L2004} gave an algorithm how
to determine a family of elliptic curves $E_t$, up to
projective equivalence.
\medskip

Yang and Yui \cite{YY2007} studied differential equations
satisfied by modular forms of two variables associated
to $\Gamma_1\times \Gamma_2$ where $\Gamma_i \, (i=1,2)$ are
genus zero subgroups of $SL(2,\RR)$ commensurable with $SL(2,\ZZ)$.
A motivation is to realize these differential
equations satisfied by modular forms of two variables 
as Picard--Fuchs differential equations
of K3 families with large Picard numbers, e.g., $19, 18, 17$ 
or $16$, thereby establishing the modularity of
solutions of Picard--Fuchs differential equations. 
This goal was achieved for some of the families of K3 surfaces studied
by Lian and Yau in \cite{LY1995}, \cite{LY1996}.
\medskip

$\bullet$ {\bf Monodromy of Picard--Fuchs differential equations
of certain families of Calabi--Yau threefolds} 

Classically, it is known that the monodromy groups of
Picard--Fuchs differential equations for families of
elliptic curves and K3 surfaces are congruence subgroups
of $SL(2,\RR)$. This modularity property of the monodromy
groups ought to be extended to families of Calabi--Yau
threefolds. For this, we will study the monodromy groups of 
Picard--Fuchs differential equations associated with one-parameter 
families of Calabi--Yau threefolds. In Chen--Yang--Yui \cite{CYY08}
they considered 14 Picard--Fuchs differential equations
of order $4$ of hypergeometric type. They are
of the form
$$\theta^4-Cz(\theta+A)(\theta+1-A)(\theta+B)(\theta+10B)$$   
where $A,B,C\in\QQ$.

\begin{theorem} {\sl
In these 14 hypergeometric cases, the matrix 
representations of the monodromy
groups relative to the Frobenius basis can be expressed in terms
of the geometric invariants of the underlying Calabi--Yau
threefolds.  Here the geometric invariants
are the {\it degree} $d$, the {\it second Chern numbers}, $c_2\cdot H$ and
the {\it Euler number},\,\, ${c_3}$. 

Furthermore, 
under suitable change of basis, the monodromy groups
are contained in certain congruence subgroups of $Sp(4,\ZZ)$
of finite index (in $Sp(4,\ZZ)$) and whose levels are related only to 
the geometric invariants.}
\end{theorem}

However, finiteness of the index of the monodromy
groups themselves in $Sp(4,\ZZ)$ is not established.

Using the same idea for the hypergeometric cases, the monodromy
groups of the differential equations of Calabi--Yau type that have
at least one conifold singularity (not of the 
hypergeometric type) are computed. Our calculations verify numerically that
if the differential equations come from geometry, then the monodromy
groups are also contained in some congruence 
subgroups of $Sp(4,\ZZ)$.

We should mention that van Enckevort and van Straten \cite{ES} numerically
determined the monodromy for $178$ Calabi--Yau equations of order $4$
with a different method from ours, and speculated that these equations
do come from geometry.
\medskip

$\bullet$ {\bf Modularity of mirror maps and mirror moonshine}

For a family of elliptic curves $y^2=x(x-1)(x-\lambda)$, the periods
$\int_1^{\infty}\frac{dx}{\sqrt{x(x-1)(x-\lambda)}}$
satisfy the Picard--Fuchs differential equation
$$(1-\lambda)\theta^2 f-\lambda\theta\,f-\frac{\lambda}{4}f=0\quad
(\theta=\lambda\frac{d}{d\lambda}).$$
The monodromy group for this Picard--Fuchs differential equation
is $\Gamma(2)\subset SL(2,\RR)$ of finite index.
The periods can be expressed in terms of the hypergeometric function
$$_2F_1(\frac{1}{2}, \frac{1}{2};1;\lambda).$$
Now suppose that $y_0(\lambda)=1+\cdots$ is the unique holomorphic
solution at $\lambda=0$ and $y_1(\lambda)=\lambda y_0(\lambda)+g(\lambda)$
be the solution with logarithmic singularity. Set 
$z=y_1(\lambda)/y_0(\lambda)$.  Then $\lambda$, as a function of $z$, 
becomes a modular function for the modular group $\Gamma(2)$. 
This is called a {\it mirror map} of the elliptic curve family.
That a mirror map is a Hauptmodul for a genus zero
subgroup $\Gamma(2)\subset SL(2,\RR)$ is referred to as
{\it mirror moonshine}.
\smallskip

For one-parameter families of K3 surfaces of generic Picard
number $19$, the Picard--Fuchs differential equations are
of order $3$.  Since such families of K3 surfaces are
equipped with Shioda--Inose structure (cf. Morrison \cite{M84}), 
the Picard--Fuchs
differential equations are symmetric squares of differential
equations of order $2$. Hence the monodromy groups are realized
as subgroups of $SL(2,\RR)$. 

Classically, explicit period and mirror maps for the quartic K3 surface
have been described by several articles, e.g., Hartmann \cite{H11} 
(and also see Lian and Yau \cite{LY96}. 
Consider the deformation of the Fermat quartic:
$F_t:=\{x_0^4+x_1^4+x_2^4+x_3^4-4tx_0x_1x_2x_3=0\}\subset\PP^3$.
Taking the quotient of $F_t$ by some finite group and then resolving
singularities, we obtain the Dwork pencil of K3 surfaces,
denoted by $X_t$. Then $X_t$ is a K3 surface
with generic Picard number $19$. 
The Picard--Fuchs differential equation is of order $3$
and there is a unique holomorphic solution (at $t=0$) of the
form $w_0(t)=1+\sum_{n\geq 1} c_nt^n$, and another solution
of logarithmic type: 
$w_1(z)=\mbox{log} w_0(t)+\sum_{n\geq 1} d_nt^n$.
Now introduce the new variable $z$ by
$z:=\frac{1}{2\pi i}\frac{w_1(t)}{w_0(t)}$, and put $q=e^{2\pi iz}$. The
inverse $t=t(q)$ is the mirror map
and is given by
$$t(q)=q-104q^2+6444q^3-311744q^4+13018830q^5+\cdots$$
This is the reciprocal of the Hauptmodul for $\Gamma_0(2)_+\subset SL(2,\RR)$. 
\smallskip

There are several more examples of one-parameter families of
K3 surfaces with generic Picard nubmer $19$. 
Doran \cite{D20} has established the modularity
of the mirror map for $M_n$-polarized K3 surfaces. However,
for each $n$, the explicit description of mirror maps as modular functions
are still to be worked out.
\medskip

The situation will get more much complicated when we
consider two-parameter families of K3 surfaces.
Hashimoto and 
Terasoma \cite{HT2011} have studied the period and mirror maps of 
the two-parameter (in fact, projective $1$-parameter) family 
$\{{\mathcal{X}}_t\}$ $t=(t_0,t_1)\in\PP^1$ of 
quartic family of K3 surfaces defined by
$${\mathcal{X}}_t: x_1+\cdots+x_5=
t_0(x_1^4+\cdots+x_5^4)+t_1(x_1^2+\cdots+x_5^2)^2=0$$
in $\PP^4$ with homogeneous coordinates 
$(x_1:\cdots :x_5)$. This family admits a symplectic group 
action by the symmetric group $S_5$. 
The Picard number of a generic fiber is 
equal to $19$, and the Gram matrix of the transcendental lattice 
$T$ is given by  
$\left(\begin{array}{ccc} 4 & 1 & 0 \\
        1 & 4 & 0 \\
        0 & 0 & -20 \end{array} \right)$.
The image of the period map of this family is a $1$-dimensional 
subdomain $\Omega_T$ of the $19$-dimensional period domain 
(the bounded symmetric domain of type IV).
Let $\Omega_T^{\circ}$ be a connected component of $\Omega_T$.
Since $O(T)$ has no cusp, there is a modular embedding
$i : \Omega_T^{\circ}\to {\bf{H}}_2$ to the Siegel upper
half-plane ${\bf{H}}_2$ of genus $2$. This modular
embedding is constructed using the Kuga-Satake
construction.
The inverse of the period map, that is, a mirror map, is
constructed using automorphic forms of one variable
on $\Omega_T^{\circ}$. 
In fact, automorphic forms are constructed as the pull-backs of
the fourth power of theta constants of genus $2$.
This gives yet another example of a generalized mirror moonshine.
\smallskip

For a one-parameter family of Calabi--Yau threefolds,
the mirror map is defined using specific solutions
of the Picard--Fuchs differential equation of the
family. At a point of maximal unipotent monodromy
(e.g., $z=0$), there is a unique holomorphic power
series solution $\omega_0(z)$ with $\omega_0(0)=1$, and a logarithmic power
series solution $\omega_1(z)=\mbox{log}(z)\omega_0(z)+g(z)$ where
$g(z)$ is holomorphic near $z=0$ with $g(0)=0$. 
Now put $t:=\frac{\omega_1(z)}{\omega_0(z)}$. We call the
map defined by $q:=e^{2\pi it}= ze^{g(z)/\omega_0(z)}$ the {\it mirror
map} of the Calabi--Yau family. (See, for instance, Lian and Yau \cite{LY96}).

For some one-parameter families of Calabi--Yau threefolds, e.g.,
of hypergeometric type, the integrality of the mirror maps
has been established by Krattenthaler and Rivoal \cite{KR2009}, 
\cite{KR2010}.

The modularity of mirror maps of Calabi--Yau families
is getting harder to deal with in general. Doran \cite{D20}
has considered certain one-parameter families
of Calabi--Yau threefolds with $h^{2,1}=1$. 
The Picard--Fuchs differential equations of
these Calabi--Yau threefolds are of order $4$. Under some
some special constraints imposed by special
geometry, and some conditions about a point $z=0$ of
maximal unipotent monodromy, there is a set of
fundamental solutions (to the Picard--Fuchs
differential equation) of the form
$\{ u, u\cdot t, u\cdot F^{\prime}, (tF^{\prime}-2F)\}$ where
$u=u(z)$ is the fundamental solution locally holomorphic at $z=0$,
$t=t(z)$ is the mirror map, $F(z)$ is the prepotential
and $F^{\prime}$ is the derivative of $F$ with respect to
$z$.  When there are no instanton corrections, the
Picard--Fuchs differential equation becomes the symmetric
cube of some second-order differential equation. In
this case, Doran has shown that the mirror map becomes 
automorphic. However, exhibiting
automorphic forms explicitly remains an open problem.
On the other hand, if there are instanton corrections,
a necessary and sufficient condition is presented for
a mirror map to be automorphic.

\medskip

We should mention here a converse approach to the modularity
question of solutions of Picard--Fuchs differential equations,
along the line of investigation by Yang and Yui \cite {YY2007}.
The starting point is modular forms and the differential equations
satisfied by them. It may happen that these differential equations 
coincide with Picard--Fuchs differential equations of some families
of Calabi--Yau varieties. Consequently, the modularity of solutions
of Picard--Fuchs differential equations and mirror maps can be
established.

\section{The modularity of generating functions
of counting some quantities on Calabi--Yau varieties}

Under this subtitle, topics included are enumerative geometry, Gromov--Witten
invariants, and various invariants counting some mathematical/physical
quantities, etc.

$\bullet$ Mirror symmetry for elliptic curves
and quasimodular forms

We consider the generating function, $F_g(q)$, counting simply 
ramified covers of genus $g\geq 1$ over a fixed elliptic curve with
$2g-2$ marked points. 

\begin{theorem} {\sl For each $g\geq 2$, $F_g(q)$ (with $q=e^{2\pi i\tau},\,
\tau\in\mathfrak{H}$), is a quasimodular form of weight $6g-6$
on $\Gamma=SL_2(\ZZ)$. Consequently, $F_g(q)$ is a polynomial
in $\QQ[E_2, E_4, E_6]$ of weight $6g-6$.}
\end{theorem}

This result is stated as the Fermion Theorem in Dijkgraaf \cite{D94},
which is concerned with the $A$-model side of mirror symmetry
for elliptic curves.
A mathematically rigorous proof was given in the article
of Roth--Yui \cite{RY2010}.
The B-model (bosonic) counting constitutes the mirror side of
the calculation. The bosonic counting will involve calculation
with Feynman integrals of trivalent graphs. A mathematical rigorous
treatment of the B-model counting is currently under way.
\smallskip

{\bf Further generalizations}: 
\smallskip

(a) The generaiting function of $m$-simple covers for any integer 
$m\geq 2$ of genus $g\geq 1$ over
a fixed elliptic curve with $2g-2$ marked points has been shown again
to be quasimodular forms by Ochiai \cite{O}.

(b) The quasimodularity of the Gromov--Witten invariants for
the three elliptic orbifolds with simple elliptic singularities
$\tilde{E}_N$ ($N=6,7,8)$ has been established by Milanov and Ruan \cite{MR}. 
These elliptic orbifolds are realized as quotients of hypersurfaces
of degree $3, 4$ and $6$ in weighted projective $2$-spaces with
weights $(1,1,1), (1,2,2)$ and $(1,2,3)$, respectively.

(c) The recent article of Rose \cite{R2012} has proved the
quasimodularity of the generating function for the number of
hyperelliptic curves (up to translation) on a polarized
abelian surface.

\section{Future Prospects}

Here we collect some topics which we are not able to cover
in this paper as well some problems for further investigation.

\subsection{The potential modularity}

The potential modularity of families of hypersurfaces.
For the Dwork families of one-parameter hypersurfaces, the
potential modularity has been established by R. Taylor and
his collaborators.  Extend the potential modularity to
more general Calabi--Yau hypersurfaces, Calabi--Yau complete
intersections, etc.

\subsection{The modularity of moduli of families of
Calabi--Yau varieties}

Moduli spaces of lattice polarized K3 surfaces
with large Picard number.

\subsection{Congruences, formal groups}

Congruences for Calabi--Yau families, formal groups.

\subsection{The Griffiths Intermediate Jacobians of Calabi--Yau
threefolds}

$\bullet$ Explicit description of the Griffiths intermediate
Jacobians and their modularity.

Let $X$ be a Calabi--Yau threefold defined over $\QQ$.
Let $$J^2(X)\simeq H^3(X,\CC)/F^2H^3(X,\CC)+H^3(X,\ZZ)\simeq
 H^3(X,\CC)^*/H_3(X,\ZZ)$$
be the Griffiths intermediate Jacobian of $X$. 
There is the Abel--Jacobi map
$$CH^2(X)_{hom,\QQ} \to J^2(X)_{\QQ}.$$
A part of the Beilinson--Bloch conjecture asserts that this
map is injective modulo torsion.

Now suppose that $X$ is rigid.
Then $J^2(X)$ is a complex torus of dimension $1$ so that
there is an elliptic curve $E$ such that $J^2(X)\simeq E(\CC)$.
We know that $X$ is modular by \cite{GY2011}. 
\medskip

{\bf Question}: {\sl Is it true that the Griffiths intermediate Jacobian
$J^2(X)$ of a rigid Calabi--Yau threefold $X$ over $\QQ$
is defined over $\QQ$ and hence modular?} 
\medskip

$\bullet$ Special values of $L$-series of Calabi--Yau
threefolds over $\QQ$.

Assuming a positive answer to the above question, we can consider
a possible relation between the Birch and Swinnerton-Dyer conjecture
for rational points on $J^2(X)_{\QQ}$: $$\mbox{rank}\,_{\ZZ} J^2(X)_{\QQ}(\QQ)
=\mbox{ord}_{s=1} L(J^2(X)_{\QQ},s)$$
and the Beilinson--Bloch conjecture on the Chow group $CH^2(X)_{\QQ}$ of $X$:
$$\mbox{rank}\,_{\ZZ} CH^2(X)_{hom,\QQ} = \mbox{ord}_{s=2} L(X_{\QQ},s).$$

If the Abel--Jacobi map $CH^2(X)_{hom,\QQ}\to H^2(X)_{\QQ}$ is injective modulo torsion, then
$$\mbox{ord}_{s=2} L(X_{\QQ},s)\leq \mbox{ord}_{s=1} L(J^2(X)_{\QQ},s).$$ 

\subsection{Geometric realization problem (The converse problem)}

We know that every singular K3 surface $X$ over $\QQ$ is motivically
modular in the sense that the transcendental
cycles $T(X)$ corresponds to a newform of weight $3$. 
For singular K3 surfaces over $\QQ$, the converse problem asks: 
\medskip

{\sl Which newform of weight $3$ with integral Fourier
coefficients would correspond to a singular
K3 surface defined over $\QQ$?} 
\medskip

This has been answered by Elkies and Sch\"utt \cite{ES2008} 
(see also the article by Sch\"utt \cite{S2012} in this volume). 
Their result is the following theorem.

\begin{theorem} {\sl Every Hecke eigenform of weight
$3$ with eigenvalues in $\ZZ$ is associated to a singular K3 surface
defined over $\QQ$.}
\end{theorem}

Now we know that every rigid Calabi--Yau threefold over
$\QQ$ is modular (see Gouv\^ea--Yui \cite{GY2011}).
The converse problem that has
been raised, independently, by Mazur and van Straten is the
so-called {\it geometric realization problem}, and is
stated as follows:
\medskip

{\sl Which newforms of weight $4$ on some $\Gamma_0(N)$ with integral
Fourier coefficients would arise from rigid
Calabi--Yau threefolds over $\QQ$? Do all such forms
arise from rigid Calabi--Yau threefolds over $\QQ$?}
\medskip

For Calabi--Yau threefolds, a very weak version of the above problem has
been addressed in Gouv\^ea--Kiming--Yui \cite{GKY2011}.
\medskip

{\bf  Question}: {\sl
Given a rigid Calabi--Yau threefold $X$ over $\QQ$ and
a newform $f$ of weight $4$, for any non-square rational
number $d$, there is a twist $f_d$ by the quadratic
character corresponding to the quadratic extension
$\QQ(\sqrt{d})/\QQ$. Does $f_d$ arise from a rigid
Calabi--Yau threefold $X_d$ over $\QQ$?} 
\medskip

A result of Gouv\^ea--Kiming--Yui \cite{GKY2011} in this
volume is that the answer is positive if a Calabi--Yau
threefold has an anti-symplectic involution. 
Let $X$ be a rigid Calabi--Yau threefold over $\QQ$. For
a square-free $d\in\QQ^{\times}$, let $K:=\QQ(\sqrt{d})$ and
let $\sigma$ be the non-trivial automorphism of $K$. We say that
a rigid Calabi--Yau threefold $X_d$ defined over $\QQ$ is
a {\it twist of $X$ by $d$} if there is an involution
$\iota$ of $X$ which acts by $-1$ on $H^3_{et}(\ov{X},\QQl)$, and
an isomorphism $\theta: (X_d)_K\cong X_K$ defined over $K$
such that $\theta^{\sigma}\circ \theta^{-1}=\iota$.

\begin{proposition} {\sl Let $X$ be a rigid Calabi--Yau threefold
over $\QQ$ and let $f$ be the newform of weight $4$ attached
to $X$.  Then, if $X_d$ is twist by $d$ of $X$,
the newform attached to $X_d$ is $f_d$, the twist of $d$ by
the Dirichlet character $\chi$ corresponding to $K$.}
\end{proposition}  

Various types of modular forms have
appeared in the physics literature.  We wish to understand
``conceptually'' why modular forms play such pivotal
roles in physics.  Here we list some of the modular
appearances in the physics literature.

\subsection{Modular forms and Gromov--Witten invariants}

Modular (automorphic) forms and Gromov--Witten invariants
and generalized invariants.  

\subsection{Automorphic black hole entropy}

This has something to do with Conformal Field Theory. Mathematically,
mock modular forms, Jacobi forms etc. will come into the picture.
This is beyond the scope of this article.

\subsection{$M_{24}$-moonshine} 

Recently, some close relations between the elliptic genus of
K3 surfaces and the Mathieu group $M_{24}$ along the line
of moonshine have been observed in the physics literature, e.g.,
\cite{EH2009a},\cite{EH2009b}. It has been observed that 
multiplicities of the non-BPS representations are given
by the sum of dimensions of irreducible representations
of $M_{24}$ and furthermore, they coincide with Fourier 
coefficients of a certain mock theta function. 
\medskip

\section{Appendix}

In this appendix, we will recall modular forms of various
kinds, e.g., classical modular forms, quasimodular forms,
Hilbert modular forms, Siegel modular forms, and most
generally, automorphic forms, which are relevant to our
discussions. For details, the reader is referred to
\cite{BGHZ08}. 

\begin{definition}
{\bf Dimension 1}:  Let $\mathfrak{H}:=\{z\in\CC\,|\, \mbox{Im}(z)>0\,\}$
be the complex upper-half plane. For a given integer $N>0$, let
$$\Gamma_0(N)=:\left\{\begin{pmatrix} a & b \\
 c & d\end{pmatrix}\in SL(2,\ZZ)\,|\, c\equiv 0\mod N\right\}$$
be a congruence subgroup of $SL(2,\ZZ)$ (of finite index).
A modular form of {\it weight} $k$ and {\it level} $N$ is a
holomorphic function $f: \mathfrak{H}\to \CC$ with the following
properties:

(M1) For $\begin{pmatrix} a & b \\
                          c & d\end{pmatrix}\in \Gamma_0(N)$,
$f\left(\frac{az+b}{cz+d}\right)=(cz+d)^kf(z)$;

(M2) $f$ is holomorphic at the cusps.

Since $\begin{pmatrix} 1 & 1 \\
                0 & 1\end{pmatrix}\in\Gamma_0(N)$, (a1)
implies that $f(z+1)=f(z)$, so $f$ has the Fourier expansion
$f(z)=f(q)=\sum_{n\geq 0} c(n)q^n$ with $q=e^{2\pi iz}$.  $f$ is a
{\it cusp} form if it vanishes at all cusps.

If $\chi$ is a mod $N$ Dirichlet character, we can define
a modular form with character $\chi$ by replacing (a1) by
$f\left(\frac{az+b}{cz+d}\right)=\chi(d)(cz+d)^kf(z)$.

The space $S_k(\Gamma_0(N))$ of all
cusp forms of weight $k$ and level $N$ is a
finite dimensional vector space, and similarly,
so also is the space $M_k(\Gamma_0(N),\chi)$ of all modular forms of
weight $k$ and level $N$.

On $S_k(\Gamma_0(N),\chi)$ there are Hecke operators $T_p$ for
every prime $p$ not dividing $N$. A cusp form $f$ is a (normalized) Hecke 
{\it eigenform} if it is an eigenvector for all $T_p$, that is,
$T_p(f)=c(p)f$. For such a normalized eigenform $f$,
define the $L$-series $L(f,s)$ by
$$L(f,s)=\sum_{n\geq 1} c(n)n^{-s}
=\prod_p \frac{1}{1-c(p)p^{-s}+\chi(p)p^{k-1-2s}}\quad\mbox{where $\chi(p)=0$
 if $p|N$}$$
\medskip

{\bf Dimension 2}: 
Let $F=\QQ(\sqrt{d})$ be a totally
quadratic field over
$\QQ$ where $d>0$ and square-free, and let $\OO_F$ be its ring of integers.
The $SL_2(\OO_F)$ can be embedded into
$SL_2(\RR)\times SL_2(\RR)$ via the two real embeddings of $F$
to $\RR$, and it acts on $\mathfrak{H}\times\mathfrak{H}$
via fractional linear transformations:
$$\begin{pmatrix} a & b \\ c & d\end{pmatrix}\,z
=\left(\frac{az_1+b}{cz_1+d}, \frac{az_2+b}{cz_2+d}\right)$$
for $z=(z_1,z_2)\in \mathfrak{H}\times\mathfrak{H}.$
The group $$\Gamma(\OO_F\oplus\mathfrak{a})\:=\left\{
\begin{pmatrix} a & b \\ c & d\end{pmatrix}\in SL_2(F), a, d\in\OO_F,
b\mathfrak{a}^{-1}, c\in\mathfrak{a}\,\right\}$$
is called {\it the Hilbert modular group} corresponding to a
fractional ideal $\mathfrak{a}$ of $F$.
If $\mathfrak{a}=\OO_F$, put
$\Gamma_F=\Gamma(\OO_F\oplus \OO_F)=SL_2(\OO_F)$.
Let $\Gamma\subset SL_2(F)$ be a subgroup commensurable with
$\Gamma_F$, and let $(k_1,k_2)\in\ZZ\times\ZZ$.

A meromorphic function $f:\mathfrak{H}\times\mathfrak{H}\to \CC$
is called a meromorphic {\it Hilbert modular form} of {\it weight}
$(k_1,k_2)$ for $\Gamma$ if
$$f(\gamma\,z)=(cz_1+d)^{k_1}(cz_2+d)^{k_2} f(z)$$
for all $\gamma=\begin{pmatrix} a & b \\ c & d\end{pmatrix}\in\Gamma$
and $z=(z_1,z_2)\in\mathfrak{H}^2$.
If $f$ is holomorphic, then $f$ is a holomorphic Hilbert modular form, and a
holomorphic Hilbert modular form is
{\it symmetric} if $f(z_1,z_2)=f(z_2,z_1)$.
Further, $f$ is a cusp form if it vanishes at cusps of $\Gamma$.
A cusp form of weight $(2,2)$ is identified with a
holomorphic $2$-form on the
Hilbert modular surface $\mathfrak{H}^2/\Gamma$.
The space of holomoprhic Hilbert modular forms of weight $(k_1,k_2)$
for $\Gamma$ is denoted by $M_k(\Gamma)$ and the cusp forms
by $S_k(\Gamma)$. $M_k(\Gamma)$ is a finite dimensional
vector space over $\CC$.

A holomorphic Hilbert modular form has a Fourier expansion at $\infty$
of the following form.
Let $M\subset F$ be a rank $2$ lattice
and $V\subset\OO^*_F$ be a finite index subgroup acting on $M$ in
a suitable way. Then
$$f(z)=a_0+\sum_{\nu\in M^{\vee},\,\nu\geq 0} 
a_{\nu}e^{2\pi i\mbox{tr}(\nu z)}$$
where $\nu$ runs over the dual lattice $M^{\vee}$, and
$\mbox{tr}(\nu\,z):=vz_1+v^{\prime} z_2$.

The $L$-series of a Hilbert modular cusp form is defined by
$$L(f,s)=\sum_{\mathfrak{a}\subset\OO_F} a(\mathfrak{a}) N(\mathfrak{a})^{-s}
$$
where $\mathfrak{a}$ runs over principal ideals.
For details, the reader should consult Bruinier \cite{Br2008}.
\medskip

{\bf Dimension 3}: 
Let $g, N\in\NN$.
Define the Siegel upper-half plane by
$$\mathfrak{H}_g:=\left\{\, z\in M_{g\times g}(\CC)\,
|\, z^t=z, \mbox{Im}(z>0\,\right\},$$
and the symplectic group $Sp(2g,\ZZ)$ as the automorphism
group of the symplectic lattice $\ZZ^{2g}$. That is,
$$Sp(2g,\ZZ)=\{\,\gamma\in GL(2g,\ZZ)\,|\, \gamma^t J\gamma=J\}$$ where
$J_g:=\begin{pmatrix} 0 & I_g \\ -I_g & 0\end{pmatrix}$. The
group $Sp(2g,\ZZ)$ acts on $\mathfrak{H}_g$ by
$$z\mapsto \gamma(z)=(Az+B)(Cz+D)^{-1}\quad\mbox{for}\quad \gamma
=\begin{pmatrix} A & B \\ C & D\end{pmatrix}\in Sp(2g,\ZZ).$$
Let
$\Gamma_g(N):=\left\{ \gamma\in Sp(2g,\ZZ)\,|\,
               \gamma\equiv I_{2g}\pmod N\,\right\}$
be a subgroup of $Sp(2g,\ZZ)$. Then $\Gamma_g(N)$ acts
freely when $N\geq 3$. Here
$I_{2g}$ is the identity matrix of order $2g$. If $N=1$,
$\Gamma_g(1)=Sp(2g,\ZZ)$.
The quotient space $\mathfrak{H}_g/\Gamma_g(N)$
(with $N\geq 3$) is a complex manifold of dimension
$g(g+1)/2$ (associated to a graded algebra of modular forms).

A holomorphic function
$f : \mathfrak{H}^g\to \CC$ is a
{\it Siegel modular form of genus $g$, weight $k\in\NN$} if
$$f(\gamma(z))=\mbox{det}(Cz+D)^kf(z)$$
for all $\gamma=\begin{pmatrix} A & B \\ C & D\end{pmatrix}\in Sp(2g,\ZZ)$ an
d all $z\in\mathfrak{H}_g$.
The space of all holomorphic Siegel modular forms
is a finitely generated graded algebra. The simplest
examples of Siegel modular forms are given by
theta constants.

A holomorphic Siegel modular form
$f$ has a Fourier expansion of the form
$$f(z)=\sum A(n)e^{2\pi \mbox{tr}(nz)}$$
where the sum runs over all positive semi-definite integral matrices
$n\in GL(g,\QQ)$. If the Fourier expansion is supported only
on positive definite integral $g\times g$ matrices $n$, then $f$ is
called a cusp form.

There are at least two different $L$-series of a holomorphic Siegel
modular form $f$: one is the spinor $L$-series, and the other is
the standard $L$-series.

However, it is not clear how to associate these $L$-series to
the Fourier expansion.

For details, the reader should consult van der Geer \cite{vG2008}.
\end{definition}

\begin{example}
(a) For $N=1$, the total space $M:=\oplus_{k} M_k(SL(2,\ZZ))$,
of all modular forms is generated by
the Eisenstein series $E_4$ and $E_6$. The Eisenstein
series $E_2$ is not modular but it may be called {\it quasimodular}.
The space $\tilde M$ of all quasimodular forms for $SL(2,\ZZ)$ is
generated by the Eisenstein series $E_2, E_4$ and $E_6$, that is,
$\tilde M=\CC[E_2, E_4, E_6]$.

For $N>1$, the space of modular forms of weight $k$ for any
congruence subgroup of level $N$ 
is finite dimensional, and its basis can be determined. 

Some properties of quasimodular forms for 
finite index subgroups of $SL_2(\ZZ)$, dimension, basis, etc.

(b) There are Eisenstein series for $\Gamma_F=SL_2(\OO_F)$ and
$k$ even given by
$$G_{k,B}=N(\mathfrak{b})^k\sum_{(c,d)\in\OO_F^*\setminus \mathfrak{b}^2}
N(cz+d)^{-k}$$
for $B$ an ideal class of $F$.

If we put $g_k:=\frac{1}{\zeta_F(k)} G_{k,\OO_F}$, then
$g_2, g_6$ and $g_{10}$ generate the graded algebra
$M_{2^*}^{symm}(\Gamma_F)$ over $\CC$, that is,
$$M_{2^*}^{symm}(\Gamma_F)\cong \CC[g_2,g_6,g_{10}].$$

(c) (Igusa) For $g=2$, the graded algebra
$\mathfrak{M}$ of classical Siegel modular forms of genus $2$
is generated by the Eisenstein series $E_4$ and $E_6$,
the Igusa cusp forms $C_{10}$, $C_{12}$,
and $C_{35}$ (where the subindex  denote weights),
and
$$\mathfrak{M}\cong
\CC[E_4,E_6,C_{10},C{12},C_{35}]/(C_{35}^2=P(E_4,E_6,C_{10},C_{12}))$$ where
$P$ is an explicit polynomial.
\end{example}

\vskip 1cm
 
\centerline{\bf Acknowledgments}

I would like to thank Matthias Sch\"utt for carefully
reading the preliminary version of this paper
and suggesting numerous improvements.
I am indebted to Ron Livn\'e for answering my  
questions about Galois representations.

I would also like to thank a number of colleagues
for their comments and suggestions. This includes
Jeng-Daw Yu, Ling Long, Ken-Ichiro Kimura, and
Bert van Geemen.

We are grateful to V. Nikulin for allowing us to
use the template of Nikulin's pryamid in
Figure 1. 

Last but not least, my sincere thanks is to the
referee for reading through the earlier versions
of this article and for very helpful constructive
criticism and suggestions for the improvement of
the article. I would also like to thank Arther Greenspoon
of Mathematical Reviews for copy-editing the article.

The article was completed while the author held
visiting professorship at various institutions
in Japan: Tsuda College, Kavli Institute for Physics
and Mathematics of the Universe and Nagoya University.
I thank the hospitality of these instutions.

\end{document}